\documentclass[12pt,a4paper,oneside]{article}
\usepackage{amsmath,amssymb}
\usepackage{amsthm}
\usepackage{accents} 
\usepackage{cite}
\usepackage{graphicx}
\usepackage[dvipsnames]{xcolor}
\usepackage{tikz}
\usepackage{subcaption}
\usepackage{caption}
\usetikzlibrary{arrows,shapes}
\usetikzlibrary{mindmap,trees}
\usetikzlibrary{shapes.misc}
\usetikzlibrary{decorations.pathmorphing, patterns}
\usetikzlibrary{decorations.pathreplacing}
\usepackage{pgfplots}
\pgfplotsset{compat=newest}
\makeatletter

\pgfdeclarepatternformonly[\LineSpace]{my north east lines}{\pgfpointorigin-0.1pt}{\pgfqpoint{\LineSpace}{\LineSpace}}{\pgfqpoint{\LineSpace}{\LineSpace}}%
{
	\pgfsetcolor{\tikz@pattern@color}
	\pgfsetlinewidth{0.4pt}
	\pgfpathmoveto{\pgfqpoint{0pt}{\LineSpace}}
	\pgfpathlineto{\pgfqpoint{100pt}{\LineSpace}}
	\pgfusepath{stroke}
}
\makeatother
\newdimen\LineSpace
\tikzset{
	line space/.code={\LineSpace=#1},
	line space=8pt
}
\numberwithin{equation}{section}

\theoremstyle{plain}
\newtheorem{theorem}{Theorem}[section]
\newtheorem{proposition}[theorem]{Proposition}
\newtheorem{definition}[theorem]{Definition}
\newtheorem{corollary}[theorem]{Corollary}

\theoremstyle{definition}
\newtheorem{remark}[theorem]{Remark}

\newcommand{\eq}{equation}
\newcommand{\real}{\ensuremath{\mathbb R}}

\newcommand{\rn}{\ensuremath{{\mathbb R}^n}}

\newcommand{\Rnp}{\ensuremath{{\mathbb R}^{n+1}}}

\newcommand{\no}{\ensuremath{\nat_0}}
\newcommand{\ganz}{\ensuremath{\mathbb Z}}
\newcommand{\zn}{\ensuremath{{\mathbb Z}^n}}

\newcommand{\Aspq}{A_{p,q}^s(\mathbb{R}^{n})}
\newcommand{\As}{\ensuremath{A^s_{p,q}}}
\newcommand{\Ao}{\ensuremath{A^{s_0}_{p,q}(\mathbb{R}^{n})}}
\newcommand{\Aso}{\ensuremath{A^{s_0}_{p,q}}}
\newcommand{\Bspq}{B_{p,q}^s(\mathbb{R}^{n})}
\newcommand{\Bs}{\ensuremath{B^s_{p,q}}}

\newcommand{\Fs}{\ensuremath{F^s_{p,q}}}

\newcommand{\nat}{\ensuremath{\mathbb N}}

\newcommand{\di}{\ensuremath{{\mathrm d}}}

\newcommand{\vp}{\ensuremath{\varphi}}
\newcommand{\hra}{\ensuremath{\hookrightarrow}}
\newcommand{\supp}{\ensuremath{\mathrm{supp \,}}}

\newcommand{\ve}{\ensuremath{\varepsilon}}
\newcommand{\vk}{\ensuremath{\varkappa}}
\newcommand{\vr}{\ensuremath{\varrho}}
\newcommand{\pa}{\ensuremath{\partial}}

\newcommand{\LA}{\ensuremath{L^r\!\As}}

\newcommand{\Lv}{\ensuremath{L_v \big( (0,T), b, X)}}

\newcommand{\wh}{\ensuremath{\widehat}}

\newcommand{\os}{\ensuremath{\overset}}

\newcommand{\cm}{\\[0.1cm]}

\newcommand{\Fou}{\mathcal{F}}

\newcommand{\np}{\frac{n}{p}}
\newcommand{\nr}{\frac{n}{r}}
\newcommand{\nh}{\frac{n}{2}}
\newcommand{\nph}{\Big(\frac{n}{p}-\frac{n}{2}\Big)_+}
\newcommand{\Tu}{\mathcal{T}_{u_0}}
\newcommand{\Lav}[3]{\ensuremath{L_{#1}((0,T),\,#2,\,#3)}}
\newcommand{\Lva}{L_{2\alpha v}\big((0,T),\,\frac{a}{2\alpha},\,A^s_{p,\,q}({\mathbb R}^n)\big)}
\newcommand{\Sm}[1]{W^{\alpha}_{#1}}

\newcommand{\eqn}{\begin{eqnarray}}
\newcommand{\eqs}{\begin{eqnarray*}}
\newcommand{\qen}{\end{eqnarray}}
\newcommand{\qes}{\end{eqnarray*}}

\providecommand{\keywords}[1]
{
	\small	
	\textbf{\textit{Keywords and Phrases.}} #1
}
\providecommand{\classification}[1]
{
	\small	
	\textbf{\textit{Math Subject Classifications.}} #1
}

\begin{document}
\title{\Large{Fractional nonlinear heat equations and characterizations of some function spaces in terms of fractional Gauss--Weierstrass semi--groups}}
\author{F. Baaske, H.-J. Schmei{\ss}er, H. Triebel}
\maketitle
\begin{center}
	\textbf{Abstract}
\end{center}	
We present a new proof of the caloric smoothing related to the fractional Gauss--Weierstrass semi--group
in Triebel-Lizorkin spaces. This property will be used to prove existence and uniqueness of mild and strong solutions
of the Cauchy problem for a fractional nonlinear heat equation. \\[1ex]
\noindent
\classification{35K30, 35K05, 46E35}\\[1ex]



\section{Introduction}   \label{S1}
The aim of the paper is twofold. First we justify the smoothing property

\begin{equation}\label{1.1}
t^{\frac{d}{2\alpha}} \| W^\alpha_t w \, \lvert A^{s+d}_{p,q} (\rn) \| \le c \, \|w \, \lvert \As (\rn) \|, \qquad 0<t \le 1, \quad d\geq 0,
\end{equation}
of the fractional Gauss--Weierstrass semi--group $W^\alpha_t$,
\begin{\eq}   \label{1.2}
	W^\alpha_t w(x) = \big( e^{-t \lvert\xi\lvert^{2\alpha}} \wh{w} \big)^\vee (x), \qquad w\in \As (\rn), \quad \alpha >0,
\end{\eq}
in Besov and Triebel-Lizorkin spaces (see Definition \ref{D2.1})
\begin{\eq}\label{1.3}
	\As (\rn), \qquad A\in \{B,F \}, \quad s\in \real \quad \text{and} \quad 1\le p,q \le \infty.
\end{\eq}
Here $\wedge$ and $\vee$ stand for the Fourier transform and its inverse, respectively.
It will be a straightforward consequence of the characterization of some of these spaces in terms of the semi--group $W^\alpha_t$.
Based on these observations we deal secondly with the Cauchy problem
\begin{align}  \label{1.4}
\pa_t u(x,t) + (-\Delta)^\alpha u(x,t) - \sum^n_{j=1} \pa_j u^2 (x,t) &=0, &&\text{$x\in \rn$, $0<t<T$}
\\  \label{1.5}
u(x,0) &= u_0 (x), &&\text{$x\in \rn$},
\end{align}
where $0<T \le \infty$ and $2\le n \in \nat$ in the  context of the semi--group $W^\alpha_t$ in \eqref{1.2} and the fractional
Laplacian
\begin{\eq}   \label{1.6}
	(-\Delta)^\alpha w = \big( \lvert\xi\lvert^{2\alpha} \wh{w} \, \big)^\vee.
\end{\eq}
Here as usual $\pa_t = \pa/\pa t$ and $\pa_j = \pa/\pa x_j$. In particular, $(-\Delta)^\alpha = (- \Delta_x )^\alpha$ refers to the
space variables. 
If $\alpha =1$ then \eqref{1.4} refers to Burgers equation. The peculiar nonlinearity  $Du^2 = \sum^n_{j=1} \pa_j u^2$ is considered as the scalar counterpart of the related
(vector--valued) non--linearity in the Navier--Stokes equations. In general the fractional Laplacian $(-\Delta)^\alpha$ is modelling dissipation (hyperdissipation if $\alpha >1$,
hyperviscousity if $\alpha =2$). In this respect the case $\alpha =\frac{n+2}{4}$ attracted special attention. We refer to \cite{KP02} and \cite{Tao09}.
The fractional Burgers equation \eqref{1.4} with $1/2 \leq\alpha<1$ has been considered in \cite{Iwa15}. For investigations of solutions of Cauchy problems for fractional 
dissipative heat equations with several types of nonlinearities we refer to \cite{MYZ08}. It turns out that both in generalized Navier-Stokes equations and in other generalized 
equations of physical and biological relevance (such as quasi-geostrophic equations, Keller-Segel equations, chemotaxis equations) the suggestion is to replace the Laplace operator 
by a (fractional) power $(-\Delta)^\alpha$ in order to achieve at adequate mathematical models.
The motivation to study the smoothing property \eqref{1.1} for fractional Gauss-Weierstrass semi--groups comes from our interest
in so--called mild solutions of \eqref{1.4}, 
\eqref{1.5} being fixed points of the operator $\Tu$,
\begin{\eq}   \label{1.7}
	\Tu u(x,t) = W^\alpha_t u_0 (x) + \int^t_0 W^\alpha_{t-\tau} Du^2 (x,\tau) \, \di \tau, \qquad x\in \rn, \quad 0<t<T
\end{\eq}
in suitable weighted Lebesgue spaces with respect to the Bochner integral $L_v \big( (0,T), b, X)$, where
$X=A_{p,q}^s(\rn),~A\in\{B,\,F\},~ s\in\real,~1\leq p\leq\infty$, $1\leq q\leq\infty,$ stands for a Besov or Triebel-Lizorkin space. This means that for some $b\in\real$ and $v$ our solution satisfies 
\begin{eqnarray}\label{Lv}
\|f\lvert\Lv\|=\left(\int_0^T t^{bv}\|f(\cdot,t)\lvert X\|^v\text{d} t\right)^{1/v} <\infty
\end{eqnarray}
if $1\leq v <\infty$ and
\eqn\label{Lvinf}
\|f\lvert\Lav{\infty}{b}{X}\|=\underset{0<t<T}{\sup}t^{b}\|f(\cdot,t)\lvert X\| <\infty
\qen
if $v=\infty$. 
Note that (after extension from $\rn\times(0,T)$ to $\Rnp$ by zero) 
\begin{\eq}
	L_v \big( (0,T), b, X) \subset S'(\Rnp) ~~~~{\rm{if}}~~~b+\frac{1}{v} < 1
\end{\eq}
(see \cite[formulae (4.17), (4.21)]{T14}) for details and an appropriate interpretation).
Here $S'(\Rnp)$ stands for the space of tempered distributions.
As far as the smoothing property \eqref{1.1} is concerned we refer to \cite[Theorem 4.1]{T14} ($\alpha =1$, classical case), \cite[Proposition 3.4]{BaS17} ($\alpha\in\nat$)
as well as \cite[Corollary 5.4]{KnS22} and \cite[Example 4.7]{KuS22} ($\alpha > 0$).
We present a different proof in Theorem \ref{T2.6} based on characterizations of $\Aspq~(s>0)$ in terms of the fractional Gauss-Weierstrass semi--group
$W_t^\alpha$ provided in Theorem \ref{T2.4}.
As far as the Cauchy problem  \eqref{1.4}, \eqref{1.5} is concerned we follow the method developed in \cite[Subsection 4.5]{T14} and \cite{Ba15} ($\alpha=1$)
as well as \cite{BaS17} and \cite{BaS19} ($\alpha\in\nat$). The main results are contained in Theorem \ref{fix2} (existence of mild and strong solutions) and
Theorem \ref{stablesol} (locally well-posedness of the Cauchy problem). Our approach allows us to deal with the above Cauchy problem for initial data $u_0$ belonging to spaces
$\Ao$ in the so-called supercritical case $s_0 > \np -2\alpha +1$, where $\alpha >1/2, 1\leq p\leq\infty,~1\leq q\leq\infty$. 
Apart from the smoothing property \eqref{1.1} a key ingredient in the proof turns out to be the mapping property of the nonlinearity $Du^2=\sum_{j=0}^n \pa_j u^2$ in \eqref{1.4}
considered in Proposition \ref{estspace} and Corollary \ref{corest} which require the condition $s>\nph$ for spatial solution spaces $\Aspq$.\\
The paper is organized as follows. Section \ref{S2} is concerned with the characterization of spaces $\Aspq$ in terms of fractional Gauss--Weierstrass semi--groups and 
the proof of the smoothing property \eqref{1.1}. The existence, uniqueness and stability of mild and strong solutions of the
Cauchy problem \eqref{1.4}, \eqref{1.5} are treated in Section \ref{S3}. In the final Section \ref{S4} we illustrate different cases of solution spaces depending on the choice of $\alpha>0$, dimension $n$ and integrability $p$. In particular we consider how close the spaces of admitted initial data approach to the so-called critical line $s_0=\frac{n}{p}-2\alpha+1$ (for more details see explanations below). Moreover, we discuss the special case when initial data $u_0$ belong to $L_p(\rn)$, $1< p <\infty$. Finally, we investigate how our results fit in the current literature and compare them with related results.

\section{Function spaces}    \label{S2}
\subsection{Definitions and basic ingredients}    \label{S2.1}
We use standard notation. Let $\nat$ be the collection of all natural numbers and $\no = \nat \cup \{0 \}$. Let $\rn$ be Euclidean $n$-space, where
$n\in \nat$. Put $\real = \real^1$.
Let $S(\rn)$ be the Schwartz space of all complex-valued rapidly decreasing infinitely differentiable functions on $\rn$ and let $S' (\rn)$ be the space of all tempered distributions on \rn.
Furthermore, $L_p (\rn)$ with $0< p \le \infty$, is the standard complex quasi-Banach space with respect to the Lebesgue measure, quasi-normed by
\begin{\eq}   \label{2.1}
	\| f \, \lvert L_p (\rn) \| = \Big( \int_{\rn} \lvert f(x)\lvert^p \, \di x \Big)^{1/p}
\end{\eq}
with the usual modification if $p=\infty$.  
Similarly $L_p (M)$ where $M$ is a Lebesgue-measurable subset of \rn. As usual, $\ganz$ is the collection of all integers and $\zn$ where $n\in \nat$ denotes the
lattice of all points $m= (m_1, \ldots, m_n) \in \rn$ with $m_k \in \ganz$. Let $Q_{j,m} = 2^{-j}m + 2^{-j} (0,1)^n$ with $j\in \ganz$ and $m\in \zn$ be the usual dyadic cubes in $\rn$,
$n \in \nat$, with sides of length $2^{-j}$ parallel to the axes of coordinates and $2^{-j} m$ as the lower left corner. 

If $\vp \in S(\rn)$ then
\begin{\eq}  \label{2.2}
	\wh{\vp} (\xi) = (\Fou \vp)(\xi) = (2\pi )^{-n/2} \int_{\rn} e^{-ix \xi} \vp (x) \, \di x, \qquad \xi \in  \rn,
\end{\eq}
denote the Fourier transform of \vp. As usual, $\Fou^{-1} \vp$ and $\vp^\vee$ stand for the inverse Fourier transform, given by the right-hand side of
\eqref{2.2} with $i$ in place of $-i$. Here $x \xi$ stands for the scalar product in \rn. Both $\Fou$ and $\Fou^{-1}$ are extended to $S'(\rn)$ in the
standard way. Let $\vp_0 \in S(\rn)$ with
\begin{\eq} \label{2.3}
	\vp_0 (x) =1 \; \text{if } \lvert x \lvert\,\le 1 \quad \text{and} \quad \vp_0 (x) =0 \; \text{if } \lvert x\lvert\, \ge 3/2.
\end{\eq}
We define the sequences
\begin{\eq}   \label{2.4}
	\vp_j (x) = \vp_0 (2^{-j} x) - \vp_0 (2^{-j+1} x ), \qquad x\in \rn, \quad j\in \nat
\end{\eq}
and
\begin{\eq}   \label{2.4a}
	\vp^j (x) = \vp_0 (2^{-j} x) - \vp_0 (2^{-j+1} x ), \qquad x\in \rn, \quad j\in \ganz.
\end{\eq}
Since
\begin{\eq}   \label{2.5}
	\sum^\infty_{j=0} \vp_j (x) =1,\;  x\in \rn \mbox{ and } \sum^\infty_{j=-\infty} \vp^j (x) =1,\;  x\in \rn\setminus\{0\},
\end{\eq}
$\vp =\{ \vp_j \}_{j\in\nat_0}$ and  $\vp =\{ \vp^j \}_{j\in\ganz}$ form a dyadic resolution of unity, respectively. The entire analytic functions $(\vp_j \wh{f} )^\vee (x)~(j\in\no)$ and $(\vp^j \wh{f} )^\vee (x)~(j\in\ganz)$ make sense pointwise in $\rn$ for any $f\in S'(\rn)$. \\
\\We are interested in inhomogeneous Besov and Triebel-Lizorkin spaces $\Aspq$ with $A\in\{B,F\}$ with $s\in\real$ and $0<p,q\leq\infty$. The standard norms of these spaces and their homogeneous counterparts are given as follows

\begin{definition}   \label{D2.1}
	Let $0<p \le \infty$ ($p<\infty$ if $A=F$), $0<q \le \infty$ and $s \in \real$.
	\\[0.1cm]
	{\em (i)} Let $\vp =\{ \vp_j \}_{j\in\nat_0}$ be the above dyadic resolution of unity.
	Then $\Aspq$ is the collection of all $f \in S' (\rn)$ such that
	\begin{\eq}   \label{2.6}
		\| f \, \lvert \Aspq \|_{\vp} = \begin{cases}
			\Big( \sum\limits^\infty_{j=0} 2^{jsq} \big\| (\vp_j \widehat{f})^\vee \, \lvert L_p (\rn) \big\|^q \Big)^{1/q},&\text{if } A=B\\
			\Big\| \Big( \sum\limits^\infty_{j=0} 2^{jsq} \big\lvert (\vp_j \wh{f})^\vee (\cdot) \big\lvert^q \Big)^{1/q} \big\lvert L_p (\rn) \Big\|,&\text{if } A=F
		\end{cases}		
	\end{\eq}
	is finite 
	$($with the usual modification if $q= \infty)$. 
	\cm
	{\em (ii)} Let $\vp =\{ \vp^j \}_{j\in\ganz}$ be the above homogeneous dyadic resolution of unity in $\rn\setminus\{0\}$.
	Then $\dot{A}^s_{p,q}(\rn)$ is the collection of all $f \in S' (\rn)$ such that
	\begin{\eq}   \label{2.7}
		\| f \, \lvert \dot{A}^s_{p,q}(\rn) \|_{\vp} = \begin{cases}
			\Big( \sum\limits^\infty_{j=-\infty} 2^{jsq} \big\| (\vp^j \widehat{f})^\vee \, \lvert L_p (\rn) \big\|^q \Big)^{1/q},&\text{if } A=B\\
			\Big\| \Big( \sum\limits^\infty_{j=-\infty} 2^{jsq} \big\lvert (\vp^j \wh{f})^\vee (\cdot) \big\lvert^q \Big)^{1/q} \big\lvert L_p (\rn) \Big\|,&\text{if } A=F
		\end{cases}	
	\end{\eq}
	is finite $($with the usual modification if $q=\infty)$.
	\cm
	{\em (iii)} Let $0<q<\infty$ and $s\in \real$. Then $F^s_{\infty, q} (\rn)$ is the  collection of all $f\in S'(\rn)$ such that
	\begin{\eq}   \label{2.8}
		\| f \, \lvert \ F^s_{\infty,q} (\rn) \|_{\vp} = \sup_{J\in \no, M\in \zn} 2^{Jn/q} \Big(\int_{Q_{J,M}} \sum_{j \ge J} 2^{jsq} \big\lvert (\vp_j \wh{f} )^\vee (x) \big\lvert^q \, \di x \Big)^{1/q}
	\end{\eq}
	with $\vp =\{ \vp_j \}_{j\in\nat_0}$ as in (i) is finite.
	\cm
	{\em (iv)} Let $0<q<\infty$ and $s\in \real$. Then $\dot{F}^s_{\infty, q} (\rn)$ is the  collection of all $f\in S'(\rn)$ such that
	\begin{\eq}   \label{2.9}
		\| f \, \lvert \ \dot{F}^s_{\infty,q} (\rn) \|_{\vp} = \sup_{J\in \ganz, M\in \zn} 2^{Jn/q} \Big(\int_{Q_{J,M}} \sum_{j \ge J} 2^{jsq} \big\lvert (\vp^j \wh{f} )^\vee (x) \big\lvert^q \, \di x \Big)^{1/q}
	\end{\eq}
	with $\vp =\{ \vp^j \}_{j\in\ganz}$ as in (ii) is finite.
\end{definition}

\begin{remark}   \label{R2.2}
	We recall that all spaces defined above are independent of the respective resolution of unity
	$\vp$ according to \eqref{2.3}--\eqref{2.5} (equivalent quasi-norms). This justifies the omission of the  subscript $\vp$ in \eqref{2.6}--\eqref{2.9} in the sequel (and any other marks in connection with equivalent quasi--norms). Note that the spaces $\Aspq$ are translation invariant. This follows easily from elementary properties of the Fourier transform and the translation invariance of $L_p$ -- spaces.
	The theory of {\em inhomogeneous} spaces, including special cases and their history may be found
	in \cite{T83}, \cite{T92}, \cite{T06} and \cite{T20}. As far as {\em homogeneous} spaces are concerned we refer to \cite[Chapter 5]{T83} 
	as well as to \cite[Definition 2.8]{T15} as far as \eqref{2.9} is concerned.
	We will use these spaces only in the context of norm equivalences. Especially for our purposes it is not necessary to discuss the usual ambiguity of {\em homogeneous} spaces. Finally, $F^s_{\infty, \infty} (\rn) =
	B^s_{\infty, \infty} (\rn)$ and $\dot{F}^s_{\infty, \infty} (\rn) =
	\dot{B}^s_{\infty, \infty} (\rn)$ as discussed  in \cite[Definition 1.1, Remark 1.2, pp.\,2--3 and p.116]{T20}.
\end{remark}

\noindent We need a few specific properties of the above defined spaces
. Let $\vp_0$ and $\vp =\{ \vp^j \}_{j\in\ganz}$ be as in \eqref{2.3} and \eqref{2.4a}. 
Let $1\le p,q \le \infty$ ($p<\infty$ for $F$--spaces) and 
$s>0$. Then
\begin{\eq}   \label{2.11}
	\begin{aligned}
		\| f \, \lvert \Bs (\rn) \| & \sim \| f\, \lvert L_p (\rn)\| +
		\|f\,\lvert\,\dot{B}^s_{p,q}(\rn)\| \\
		& \sim \| (\vp_0 \wh{f} )^\vee \lvert L_p (\rn) \| +\|f\,\lvert\,\dot{B}^s_{p,q}(\rn)\| 
	\end{aligned}
\end{\eq}
are equivalent norms in $\Bs (\rn)$ and
\begin{\eq}   \label{2.12}
	\begin{aligned}
		\| f \, \lvert \Fs (\rn) \| & \sim \| f \, \lvert L_p (\rn) \|  +\|f\,\lvert\,\dot{F}^s_{p,q}(\rn)\|
		\\
		& \sim \| (\vp_0 \wh{f} )^\vee \lvert  L_p (\rn) \| + \|f\,\lvert\,\dot{F}^s_{p,q}(\rn)\|
	\end{aligned}
\end{\eq}
are equivalent norms in $\Fs (\rn)$ (with the usual modification if $q=\infty$). This is a special case of corresponding assertions
in \cite[Section 2.3.3, pp.\,97--100]{T92} where one finds also continuous versions with $t>0$ in place of $2^{-j}$, $j\in \ganz$,
which are nearer to what follows. These norms are characterizing what means that $f\in S'(\rn)$ belongs to $\As (\rn)$ if, and only
if, the corresponding norm is finite. We need an extension of the above norms to a wider class of functions $\vp^j (x) = \vp^0 (2^{-j}x)
$ and $\vp^0 (tx)$.\\
\\
Let $h\in S(\rn)$ and $H\in S(\rn)$ with
\begin{\eq}   \label{2.13}
	h(x) =1 \ \text{if $\lvert x\lvert\, \le 1$}, \qquad  \supp h \subset \{ x: \, \lvert x\lvert\, \le 2 \}
\end{\eq}
and
\begin{\eq}   \label{2.14}
	H(x) =1 \ \text{if $1/2 \le \lvert x\lvert \, \le 2$}, \qquad  \supp H \subset \{ x: \, 1/4 \le \lvert x\lvert\, \le 4 \}.
\end{\eq}  

\begin{proposition}    \label{P2.3}
	Let $\vp_0$ be as in \eqref{2.3} and $\vp \in C^\infty (\rn \setminus \{0 \} )$ with $\lvert\vp (x)\lvert >0$ if $1/2 \le \lvert x\lvert\, \le 2$. 
	\cm
	{\em (i)} Let $1 \le p,q \le \infty$ and $0<s<\sigma$. Let
	\begin{\eq}   \label{2.15}
		\int_{\rn} \Big\lvert \Big( \frac{\vp (z) \, h(z)}{\lvert z\lvert^\sigma} \Big)^\vee (y) \Big\lvert \, \di y <\infty
	\end{\eq}
	and
	\begin{\eq}   \label{2.16}
		\sup_{m\in \nat} \int_{\rn} \Big\lvert \Big( \vp \big( 2^m \cdot \big) H(\cdot) \Big)^\vee (y) \Big\lvert \, \di y <\infty.
	\end{\eq}
	Then
	\begin{\eq}   \label{2.17}
		\begin{aligned}
			\|f \,\lvert\Bs (\rn)\| & \sim \|f \, \lvert L_p (\rn) \|
			+ \Big( \int^\infty_0 t^{-sq} \big\| \big( \vp (t\cdot) \wh{f} \big)^\vee \lvert L_p (\rn) \big\|^q \frac{\di t}{t} \Big)^{1/q} \\
			&\sim \big\| (\vp_0 \wh{f} )^\vee \lvert L_p (\rn) \big\|
			+ \Big( \int^\infty_0 t^{-sq} \big\| \big( \vp (t\cdot) \wh{f} \big)^\vee \lvert L_p (\rn) \big\|^q \frac{\di t}{t} \Big)^{1/q} 
		\end{aligned}
	\end{\eq}
	$($usual modification if $q= \infty)$ are equivalent norms in $\Bs (\rn)$.
	\cm
	{\em (ii)} Let $1 \le p<\infty$, $1 \le q\le \infty$, $0<s<\sigma$ and $a>n$. Let
	\begin{\eq}   \label{2.18}
		\int_{\rn} \Big\lvert \Big( \frac{\vp (z) \, h(z)}{\lvert z\lvert^\sigma} \Big)^\vee (y) \Big\lvert \, (1+\lvert y\lvert)^{a} \, \di y <\infty
	\end{\eq}
	and
	\begin{\eq}    \label{2.19}
		\sup_{m\in \nat} \int_{\rn} \Big\lvert \Big( \vp (2^m \cdot) \, H(\cdot) \Big)^\vee (y) \Big\lvert \, (1 + \lvert y\lvert)^{a} \, \di y <\infty.
	\end{\eq}
	Then
	\begin{\eq}   \label{2.20}
		\begin{aligned}
			& \|f \,\lvert\Fs (\rn)\| \\ & \sim \|f \, \lvert L_p (\rn) \|
			+\Big\| \Big( \int^\infty_0 t^{-sq} \big\lvert \big( \vp (t\cdot) \wh{f} \big)^\vee (\cdot) \big\lvert^q \frac{\di t}{t} \Big)^{1/q}
			\lvert L_p (\rn) \Big\| \\
			&\sim \big\| (\vp_0 \wh{f} )^\vee \lvert L_p (\rn) \big\|
			+\Big\| \Big( \int^\infty_0 t^{-sq} \big\lvert \big( \vp (t\cdot) \wh{f} \big)^\vee (\cdot) \big\lvert^q \frac{\di t}{t} \Big)^{1/q}
			\lvert L_p (\rn) \Big\| \\
		\end{aligned}
	\end{\eq}
	$($usual modification if $q= \infty)$ are equivalent norms in $\Fs (\rn)$.
\end{proposition}

\begin{proof}
	The extension of \eqref{2.11}, \eqref{2.12} from $\vp^0 (2^{-j}x)$ and its continuous counterpart $\vp^0 (tx)$ to the above assertion
	follows from \cite[Proposition 2.10, pp.\,18--19]{T15} and the references  given there specified to the above values of the parameters
	$s,p,q$. Again these norms are characterizations as explained above.
\end{proof}

\begin{remark}\label{normequ}
	Note that by \cite[Proposition 2.10]{T15} the second summands on right-hand-sides in \eqref{2.17} and \eqref{2.20} are equivalent norms in $\dot{B}^s_{p,q}(\rn)$ and $\dot{F}^s_{p,q}(\rn)$, respectively, for all admitted parameters.
\end{remark}

Let us shortly comment on the conditions with respect to $\vp_0$ and $\vp$ supposed in Proposition \ref{P2.3}.
The system $(\vp_j)_{j=1}^\infty$ introduced in order to define the spaces $\Aspq$ (see Definition \ref{D2.1}) can be rewritten as 
\[
\vp_j(x) = \vr (2^{-j} x),~~~\mbox{where}~~~\vr(x) = \vp_0(x) - \vp_0 (2x)
\]
and $\vp_0$ has the meaning of \eqref{2.3}. The function $\vr$ has compact support in $\{ x:~ \frac{1}{2} \leq \lvert x\lvert\, \leq \frac{3}{2}\}$ and satisfies the 
Tauberian condition $\lvert\vr(x)\lvert > 0$ on $\{ x:~ \frac{3}{4} \leq \lvert x\lvert \,\leq 1\}$. The characterization of spaces $\Aspq$ in Proposition \ref{P2.3} can be considered as
a continuous extension and generalization of Definition \ref{D2.1}, where the generating function $\vr$ is replaced by $\vp$. In contrast to the properties of $\vr$
it is not assumed that $\vp$ has compact support in a subset of $\rn \setminus \{0\}$. Conditions \eqref{2.15} and \eqref{2.18} ensure sufficiently strong decay to $0$ near the origin,
whereas \eqref{2.16} and \eqref{2.19} are responsible for decay if $\lvert x\lvert \to\infty$.
For example, it follows from \eqref{2.15} and \eqref{2.18} that 
$\lvert \vp (x)\lvert~\lesssim ~ \lvert x\lvert ^\sigma$ in a neighbourhood of the origin. Moreover, the condition $\lvert \vp (x)\lvert >0$ if $\frac{1}{2}\leq \lvert x\lvert\,\leq 2$ corresponds 
to the Tauberian condition with respect to $\vr$. Let us also mention that the condition with respect to $\vp_0\in C_0^\infty (\rn)$ can be weakened. 
For a more detailed discussion we refer to \cite[Corollary 2.4.1/1, Remark 2.4.1/3]{T92}.
Relevant examples will be treated in the next subsection.

\subsection{Characterizations of some function spaces in terms of fractional Gauss--Weierstrass semi--groups}   \label{S2.2}
We wish to apply Proposition \ref{P2.3} to
\begin{\eq}  \label{2.21}
	\vp (\xi) = \lvert\xi\lvert^\delta e^{-\lvert\xi\lvert^{2\alpha}}, \qquad \xi \in \real, \quad \alpha >0, \quad \delta >0.
\end{\eq}
If $0<\alpha \not\in \nat$ then $\vp(\xi)$ is not smooth at $\xi =0$ and some extra care is needed. This is just the reason why we
prefer now Proposition \ref{P2.3} (under the indicated restrictions for the underlying spaces) compared with the original inhomogeneous
versions according to \cite[Theorems 2.4.1, 2.5.1, pp.\,100, 101, 132]{T92} (which apply to all spaces $\As (\rn)$ with exception of
$F^s_{\infty,q} (\rn)$). Rescue comes from the following observations in \cite{MYZ08}. Let
\begin{\eq}   \label{2.22}
	K^\alpha (x) = \big( e^{-\lvert\xi\lvert^{2\alpha}} \big)^\vee (x), \qquad x\in \rn, \quad \alpha >0,
\end{\eq}
and according to \eqref{1.6}
\begin{\eq}   \label{2.23}
	K^{\alpha, \sigma} (x) = (-\Delta)^{\sigma/2} K^\alpha (x) = \big(\lvert\xi\lvert^\sigma e^{-\lvert\xi\lvert^{2 \alpha}} \big)^\vee (x), \qquad x\in \rn,
	\quad \sigma >0, \quad \alpha >0.
\end{\eq}
Then the estimates
\begin{\eq}     \label{2.24}
	\lvert K^\alpha (x) \lvert \le c \, (1+\lvert x\lvert)^{-n-2\alpha}, \qquad x\in \rn, \quad \alpha >0,
\end{\eq}
and
\begin{\eq}   \label{2.25}
	\lvert K^{\alpha, \sigma} (x) \lvert \le c \, (1+\lvert x\lvert)^{-n-\sigma}, \qquad x\in \rn, \quad \alpha >0, \quad \sigma >0,
\end{\eq}
are covered by \cite[Lemma 2.1, Lemma 2.2, pp.\,463, 465]{MYZ08}. With $W^\alpha_t$ as in \eqref{1.2} one has for $k\in \no$,
\begin{\eq}   \label{2.26}
	\pa^k_t W^\alpha_t w(x) = (-1)^k \big( \lvert\xi\lvert^{2k\alpha} \, e^{-t \lvert\xi\lvert^{2\alpha}} \wh{w} \big)^\vee (x), \qquad x\in \rn, \quad t>0.
\end{\eq}
In the distinguished case $\alpha =1$ one has now final characterizations of all spaces $\As (\rn)$ with $s\in \real$ and $0<p,q \le
\infty$ in terms of $\pa^k_t W_t w$ for the classical Gauss--Weierstrass semi--group $W_t = W^1_t$. This may be found in \cite[Section
3.2.7, pp.\,106--109]{T20} and the references given there. We extend now these assertions to the fractional Gauss--Weierstrass
semi--group $W^\alpha_t$ under the same restrictions for the spaces $\As (\rn)$ as in Proposition \ref{P2.3}.

\begin{theorem}    \label{T2.4}
	Let $\vp_0$ be as in \eqref{2.3} and
	$W^\alpha_t$ be as in \eqref{1.2} with $\alpha >0$.
	\cm
	{\em (i)} Let $1\le p,q \le \infty$, $s>0$ and $k\in \nat$ such that $2\alpha k>s$. Then
	\begin{\eq}   \label{2.27}
		\begin{aligned}
			\|f \, \lvert \Bs (\rn) \| &\sim \|f \, \lvert L_p (\rn) \| + \Big( \int^\infty_0 t^{(k - \frac{s}{2\alpha})q} \big\| \pa^k_t W^\alpha_t f \, \lvert
			L_p (\rn) \big\|^q \, \frac{\di t}{t} \Big)	^{1/q} \\
			&\sim \| (\vp_0 \wh{f}\, )^\vee \lvert L_p (\rn)\| + \Big( \int^\infty_0 t^{(k - \frac{s}{2\alpha})q} \big\| \pa^k_t W^\alpha_t f \, \lvert
			L_p (\rn) \big\|^q \, \frac{\di t}{t} \Big)	^{1/q} 
		\end{aligned}
	\end{\eq}
	$($equivalent norms$)$, usual modification if $q=\infty$.
	\cm
	{\em (ii)} Let $1\le p <\infty$, $1\le q \le \infty$, $s>0$ and $k\in \nat$ such that $2\alpha k > s+n$. Then
	\begin{\eq}   \label{2.28}
		\begin{aligned}
			&\|f \, \lvert \Fs (\rn) \| \\ &\sim \|f \, \lvert L_p (\rn) \| + \Big\|
			\Big( \int^\infty_0 t^{(k - \frac{s}{2\alpha})q} \big\lvert \pa^k_t W^\alpha_t f (\cdot) \big\lvert^q \frac{\di t}{t} \Big)^{1/q} \lvert
			L_p (\rn) \Big\| \\
			&\sim \| (\vp_0 \wh{f}\, )^\vee \lvert L_p (\rn)\| +  \Big\|
			\Big( \int^\infty_0 t^{(k - \frac{s}{2\alpha})q} \big\lvert \pa^k_t W^\alpha_t f (\cdot) \big\lvert^q \frac{\di t}{t} \Big)^{1/q} \lvert
			L_p (\rn) \Big\| 
		\end{aligned}
	\end{\eq}
	$($equivalent norms$)$, usual modification if $q=\infty$.
\end{theorem}

\begin{proof}
	{\em Step 1.} We rely on part (i) of Proposition \ref{P2.3} choosing there
	\begin{\eq}   \label{2.29}
		\vp (\xi) = \lvert\xi\lvert^\sigma e^{-\lvert\xi\lvert^{2\alpha}}, \qquad \xi \in \rn, \quad \sigma =2\alpha k >s.
	\end{\eq}
	Then \eqref{2.15} follows from \eqref{2.24} and
	\begin{\eq}    \label{2.30}
		\int_{\rn} \Big\lvert \big( e^{-\lvert\cdot\lvert^{2\alpha}} h (\cdot) \big)^\vee (y) \Big\lvert \, \di y \le \int_{\rn} \Big\lvert \big( e^{-\lvert\cdot\lvert^{2\alpha}}\big)^\vee (y) \Big\lvert \, \di y \cdot \int_{\rn} \lvert h^\vee (y)\lvert \, \di y <\infty.
	\end{\eq}
	Secondly we have to justify \eqref{2.16} with $\vp$ as in \eqref{2.29}. But this follows from
	\begin{\eq}   \label{2.31}
		\int_{\rn} \lvert g^\vee (x)\lvert \, \di x \le c \Big(\int_{\rn} \big\lvert (1+\lvert x\lvert^2)^{l/2} g^\vee (x) \big\lvert^2 \, \di x \Big)^{1/2} \sim \| g \, \lvert
		W^l_2 (\rn) \|, 
	\end{\eq}
	$n/2 <l \in \nat$, where $W^l_2 (\rn)$ are the classical Sobolev spaces. Then the second terms in \eqref{2.17} are equivalent to
	\begin{\eq}   \label{2.32}
		\begin{aligned}
			&\sim \Big(\int^\infty_0 \tau^{-sq} \big\| \big(\vp(\tau \cdot) \wh{f} \big)^\vee \lvert L_p (\rn) 
			\big\|^q \frac{\di \tau}{\tau} \Big)^{1/q} \\
			&\sim \Big( \int^\infty_0 t^{- \frac{s}{2\alpha}q} \big\|  \big( \vp (t^{\frac{1}{2\alpha}} \cdot ) \wh{f} \big)^\vee \lvert L_p (\rn)
			\big\|^q \, \frac{\di t}{t} \Big)^{1/q}
		\end{aligned}
	\end{\eq}
	again with $\vp$ as in \eqref{2.29} and $\tau = t^{\frac{1}{2\alpha}}$. One has by 
	\begin{\eq}   \label{2.33}
		\vp \big( t^{\frac{1}{2\alpha}} \xi \big) = t^k \lvert\xi\lvert^{2k\alpha} e^{-t \lvert\xi\lvert^{2\alpha}}, \quad t>0, \quad \xi \in \rn,
	\end{\eq}
	and \eqref{1.2} that
	\begin{\eq}   \label{2.34}
		\big( \vp (t^{\frac{1}{2\alpha}} \cdot) \wh{f} \big)^\vee (x) = (-1)^k t^k \pa^k_t W^\alpha_t f(x).
	\end{\eq}
	Inserted in \eqref{2.32} one obtains \eqref{2.27}.
	\cm
	{\em Step 2.} For the proof of part (ii) we rely on part (ii) of Proposition \ref{P2.3} choosing 
	\begin{\eq}   \label{2.35}
		\vp (\xi) = \lvert\xi\lvert^\delta e^{-\lvert\xi\lvert^{2\alpha}}, \qquad \xi \in \rn, \quad \delta =2\alpha k >s +n.
	\end{\eq}
	Using $(1+\lvert y\lvert)^{a} \le (1+ \lvert y-z\lvert)^{a} (1+\lvert z\lvert)^{a}$, $a>n$, one obtains similarly as in \eqref{2.30} that the expression \eqref{2.18}
	can be estimated from above by
	\begin{\eq}   \label{2.36}
		c \int_{\rn} \Big\lvert \Big( \frac{\vp (z)}{\lvert z\lvert^\sigma} \Big)^\vee (1+\lvert y\lvert)^{a} \, \di y 
		= c \int_{\rn} \big\lvert \big( \lvert z\lvert^{\delta-\sigma} e^{-\lvert z\lvert^{2\alpha}} \big)^\vee (y) \big\lvert (1+\lvert y\lvert)^{a} \, \di y
	\end{\eq}
	with $\sigma >s$ such that also $\delta - \sigma >a>n$ and some $c>0$. Now \eqref{2.18} follows from \eqref{2.23}, \eqref{2.25}. As far 
	as the terms \eqref{2.19} are concerned one argues as in \eqref{2.31} incorporating the factor $(1+\lvert x\lvert)^{a}$. Afterwards one is in 
	the same position as in \eqref{2.32}--\eqref{2.34} now based on \eqref{2.20}. This proves \eqref{2.28}.
\end{proof}

\begin{remark}  \label{R2.5}
	As already mentioned above the equivalent norms in \eqref{2.27} and \eqref{2.28} are characterizations. This means in our case that
	$\As (\rn)$ collects all $f\in L_p (\rn)$, $1\le p \le \infty$, such that the corresponding norm is finite. In particular it follows
	from the above considerations immediately that always $\pa^k_t W^\alpha_t f \in L_p (\rn)$ if $f\in L_p (\rn)$, $1\le p \le \infty$.
	But we will not stress this point in the sequel.
\end{remark}

\subsection{Smoothing properties}   \label{S2.3}
We justify \eqref{1.1}--\eqref{1.3}. As already mentioned in the Introduction assertions of this type are not new.  A proof of
\eqref{1.1} with $\alpha =1$ for the classical Gauss--Weierstrass semi--group $W_t w = W^1_t w$ covering all spaces $\As (\rn)$,
$A \in \{B,F \}$, $s\in \real$ and $0<p,q \le \infty$ may be found in \cite[Theorem 3.35, p.\,110]{T20}. It relies on 
characterization  of these spaces in terms of $W_t w$ using in a decisive way that the underlying kernel $e^{-\lvert\xi\lvert^2} \in S(\rn)$
is smooth at the origin $\xi =0$. This is no longer the case in general if one steps from $W_t$ to $W^\alpha_t$, $\alpha >0$, this
means from $e^{-\lvert\xi\lvert^2}$ to $e^{-\lvert\xi\lvert^{2\alpha}}$. On the other hand, \eqref{1.1} for the classical  Gauss--Weierstrass semi--group
$W_t = W^1_t$, restricted to $A \in \{B,F \}$, $s\in \real$ and $1\le p,q \le \infty$ ($p<\infty$ for $F$--spaces) is also a special
case of a corresponding assertion for related hybrid spaces $\LA (\rn)$. This may be found in \cite[Theorem 4.1, p.\,114]{T14}
including related references and comments. The extension of \eqref{1.1} from $\alpha =1$ to $\alpha \in \nat$ for the spaces $A\in 
\{B,F \}$, $s\in\real$ and $1\le p,q \le \infty$ ($p<\infty$ for $F$--spaces) goes back to \cite[Theorem 3.5, p.\,2123]{BaS17}. The
arguments both in \cite{T14} (including underlying references) and \cite{BaS17} rely on the elaborated machinery of (caloric) wavelet
expansions. The step from $\alpha \in \nat$ to $\alpha >0$ in \eqref{1.1} for $A\in \{B,F \}$, $s\in \real$ and $1\le p,q \le \infty$
is covered by the recent paper \cite{KuS22} in the larger context of convolution inequalities in these spaces. What follows may be
considered as a surprising simple proof of these assertions relying on Theorem \ref{T2.4} and a few well--known properties of the
spaces $\As (\rn)$ as introduced in Definition \ref{D2.1}.

\begin{theorem}   \label{T2.6}
	Let $W^\alpha_t$ be as in \eqref{1.2}. Let $A \in \{B,F \}$, $s\in \real$ and $1\le p,q \le \infty$. Let $d \geq 0$. Then there is a
	constant $c>0$ such that for all $t$ with $0<t \le 1$ and all $w\in \As (\rn)$,
	\begin{\eq}  \label{2.37}
		t^{\frac{d}{2\alpha}} \, \| W^\alpha_t w\, \lvert A^{s+d}_{p,q} (\rn) \| \le c \| w\, \lvert \As (\rn) \|.
	\end{\eq}
\end{theorem}

\begin{proof}
	{\em Step 1.} 
	Let $s>0$ and let $w\in \Bspq \subset L_p(\rn)$. We put
	\begin{\eq}   \label{2.38}
		\|w \, \lvert\accentset{\ast}{B}^s_{p,q} (\rn) \| = \Big( \int^\infty_0 t^{(k - \frac{s}{2\alpha})q} \big\| \pa^k_t W^\alpha_t w \, \lvert
		L_p (\rn) \big\|^q \, \frac{\di t}{t} \Big)	^{1/q} 
	\end{\eq}
	for the second summand on the right--hand side of \eqref{2.27}.
	According to \eqref{1.2} and \eqref{2.24} we have $f= W^\alpha_\tau w\in L_p(\rn)$ (see also Remark \ref{R2.5}) and it holds
	\begin{\eq}   \label{2.39}
		\pa^k_t W^\alpha_t f (x) = \pa^k_t\Big(e^{-t\lvert\xi\lvert^{2\alpha}}\widehat{f}\Big)^\vee(x)=
		(-1)^k \big( \lvert\xi\lvert^{2\alpha k} e^{-t\lvert\xi\lvert^{2\alpha}} \wh{f} \big)^\vee (x)~.
	\end{\eq}
	Inserting 
	\[
	f=W^\alpha_\tau w = \big( e^{-\tau\lvert\xi\lvert^{2\alpha}}\wh{w} \big)^\vee  
	\]
	we get
	\begin{\eq}   \label{2.39a}
		\begin{aligned}
			&\pa^k_t W^\alpha_t (W^\alpha_\tau w)(x) = 	
			(-1)^k \big( \lvert\xi\lvert^{2\alpha k} e^{-t\lvert\xi\lvert^{2\alpha}} e^{-\tau\lvert\xi\lvert^{2\alpha}}\wh{w}(\xi) \big)^\vee (x)\\
			& = (-1)^k \big(  \lvert\xi\lvert^{2\alpha k} e^{-(t+\tau)\lvert\xi\lvert^{2\alpha}} \wh{w}(\xi) \big)^\vee (x)
			= \big( \pa_t^k e^{-(t+\tau)\lvert\xi\lvert^{2\alpha}} \wh{w}(\xi)\big)^\vee (x)\\
			& = \pa^k_t W^\alpha_{t+\tau} w (x)~.
		\end{aligned}
	\end{\eq}
	Note that \eqref{2.39} is well defined due to \eqref{2.24}. Combining \eqref{2.38} and \eqref{2.39a} one obtains
	
	\begin{\eq}   \label{2.40}
		\tau^{\frac{d}{2\alpha}} \|W^\alpha_\tau w \, \lvert \accentset{\ast}{B}^{s+d}_{p,q} (\rn)\| 
		= \Big( \int^\infty_0 \tau^{\frac{d}{2\alpha}q} t^{(k- \frac{s+d}{2\alpha})q} \big\| \pa^k_t W^\alpha_{t+\tau} w \, \lvert L_p (\rn) 
		\big\|^q \, \frac{\di t}{t} \Big)^{1/q}.
	\end{\eq}
	Let $d\geq 0$ and let $ \frac{s+d}{2\alpha} + \frac{1}{q} <k \in \nat$. Then $a=k - \frac{s}{2\alpha} - \frac{1}{q} > \frac{d}{2\alpha}$,
	\begin{\eq}   \label{2.41}
		0 \leq \vk = \frac{d}{2 \alpha a} <1 \quad \text{and} \quad \big( k - \frac{s+d}{2\alpha} - \frac{1}{q} \big) \frac{1}{a} =1 - \vk.
	\end{\eq}
	Then it follows from $\tau^{\vk} t^{1-\vk} \le \tau +t$ that for $1\le q <\infty$
	\begin{\eq}   \label{2.42}
		\tau^{\frac{d}{2\alpha} q} \, t^{(k- \frac{s+d}{2\alpha})q -1} \le (\tau +t)^{(k- \frac{s}{2\alpha})q-1}
	\end{\eq}
	(modification if $q=\infty$). Inserted in \eqref{2.40} one obtains
	\begin{\eq}   \label{2.43}
		\tau^{\frac{d}{2\alpha}} \|W^\alpha_\tau w\, \lvert \accentset{\ast}{B}^{s+d}_{p,q} (\rn) \| \le \| w\, \lvert \accentset{\ast}{B}^s_{p,q} (\rn) \|.
	\end{\eq}
	As far as the first terms on the right--hand side of \eqref{2.27} are concerned it is sufficient to justify
	\begin{\eq}   \label{2.44}
		\big\| \big( e^{-\tau \lvert\xi\lvert^{2\alpha}} \wh{w} \big)^\vee \lvert L_p (\rn) \| \le c \, \|w \, \lvert L_p (\rn) \|
	\end{\eq}
	for some $c>0$ and all $0<\tau \le 1$. Recall that $1\le p \le \infty$. Then
	\eqref{2.44} follows from
	\begin{\eq}   \label{2.45}
		\int_{\rn} \Big\lvert \big( e^{-\lvert\lambda \xi\lvert^{2\alpha}} \big)^\vee(x) \Big\lvert \, \di x \le C
	\end{\eq}
	for some $C>0$ and all $0<\lambda <\infty$ what in turn can be obtained from \eqref{2.22}, \eqref{2.24} and
	\begin{\eq}  \label{2.46}
		\int_{\rn} \Big\lvert  \big( e^{-\lvert\lambda \xi\lvert^{2\alpha}} \big)^\vee(x) \Big\lvert \, \di x = \lambda^{-n} \int_{\rn} \Big\lvert  
		\big(e^{-\lvert\xi\lvert^{2\alpha}} \big)^\vee (\lambda^{-1} x) \Big\lvert \, \di x.
	\end{\eq}
	Now \eqref{2.37} can be obtained  for $\Bs (\rn)$ with $s>0$ and $1\le p,q \le \infty$ from \eqref{2.27}, \eqref{2.43} and \eqref{2.44}.\\
	Next we consider the case of $F$--spaces. Let $s>0$ and let $\omega\in\Fs(\rn)\subset L_p (\rn)$, where $1\leq p<\infty$. We put
	\begin{\eq} \label{2.38a}
		\|w \, \lvert\accentset{\ast}{F}^s_{p,q} (\rn) \| =\Big\| \Big( \int^\infty_0 t^{(k - \frac{s}{2\alpha})q} \,\Big\lvert \pa^k_t W^\alpha_t w(\cdot)\Big\lvert ^q\, \frac{\di t}{t} \Big)^{1/q} \, \lvert L_p (\rn)\Big\|
	\end{\eq}	
	Again it holds \eqref{2.39a}. The counterpart of \eqref{2.40} reads as
	\begin{\eq}   \label{2.41a}
		\tau^{\frac{d}{2\alpha}}\| W^\alpha_t w\,\lvert \accentset{\ast}{F}^{s+d}_{p,q} (\rn) \|=\Big\|\Big(\int_0^\infty \tau^{\frac{d}{2\alpha}q}t^{(k-\frac{s+d}{2\alpha})q}\, \Big\lvert \pa^k_t W^\alpha_t w(\cdot)\Big\lvert ^q \frac{\di t}{t} \Big)^{1/q} \, \Big\lvert L_p (\rn)\Big\|,
	\end{\eq}	
	where $d\geq 0$ and $ \frac{s+d}{2\alpha} + \frac{1}{q} <k \in \nat$. By the same arguments as in the proof of \eqref{2.43} one obtains
	\begin{\eq}   \label{2.43a}
		\tau^{\frac{d}{2\alpha}}\| W^\alpha_t w\,\lvert \accentset{\ast}{F}^{s+d}_{p,q} (\rn) \|\leq \|w\,\lvert \accentset{\ast}{F}^{s}_{p,q} (\rn) \|.
	\end{\eq}	
	Now \eqref{2.37} for $A=F$ is a consequence of \eqref{2.28}, \eqref{2.44} and \eqref{2.43a}.
	\cm
	{\em Step 2.} Recall
	\begin{\eq}   \label{2.47}
		I_\sigma \As (\rn) = A^{s+\sigma}_{p,q} (\rn), \qquad s\in \real, \quad \sigma \in \real \quad \text{and} \quad 0<p,q \le \infty,
	\end{\eq}
	$A \in \{B,F \}$, where
	\begin{\eq}    \label{2.48}
		I_\sigma f = \big( w_{-\sigma} \wh{f} \big)^\vee, \qquad f\in S'(\rn),
	\end{\eq}
	is the well--known lift based on $w_\delta (x) = (1+\lvert x\lvert^2)^{\delta/2}$, $x\in \rn$, $\delta \in \real$, 
	\cite[Section 1.3.2, p.\,16]{T20} and the references given there. 
	By definition of $I_\sigma$ and $W_t^\alpha$ it is not difficult to see that
	\begin{\eq}   \label{2.49}
		W^\alpha_t f = I_{-\sigma} \big(W^\alpha_t (I_\sigma f)\big)
	\end{\eq}
	if $f\in \Aspq,~1\leq p<\infty, ~s>0$ and $\sigma >0 $ (see Remark \ref{lift} below). If $s\leq 0$ and $f\in\Aspq$ then we take \eqref{2.49} as definition of $W^\alpha_t f$,
	where $\sigma $ is chosen such that $s+\sigma >0$. Then one can extend \eqref{2.37} from the spaces $\As (\rn)$, $s>0$ 
	and $1\le p,q \le \infty$ ($p<\infty$ for $F$--spaces) treated in
	Step 1 to their counterparts with $s\leq 0$.
	This covers all spaces in the  above theorem with exception of $F^s_{\infty,q} (\rn)$, $1\le q <\infty$.
	\cm
	{\em Step 3.} Let $\os{\circ}{F}{}^s_{p,q} (\rn)$, $1\le q \le \infty$, be the completion of $S(\rn)$ in $\Fs (\rn)$. Then one has
	\begin{\eq}   \label{2.50}
		\os{\circ}{F}{}^s_{1,q} (\rn)' = F^{-s}_{\infty, q'} (\rn), \qquad 1\le q \le \infty, \quad \frac{1}{q} + \frac{1}{q'} =1 \quad 
		\text{and} \quad s\in \real,
	\end{\eq}
	for the related dual spaces in the framework of the dual pairing $\big( S(\rn), S' (\rn) \big)$. This is a special case of \cite[
	(1.25), p.\,5]{T20} with a reference to \cite[Theorem 4, p.\,87]{Mar87} as far as the case $q=\infty$ is concerned (if $q<\infty$
	then $S(\rn)$ is already dense in $F^s_{1,q} (\rn)$ and the related duality is well known, \cite[p.\,5]{T20} and the references
	there). If $\vp \in S(\rn)$ then $W^\alpha_t \vp$ can be approximated in, say, $F^{s+1}_{1,1} (\rn)$ by functions belonging to 
	$S(\rn)$ for any $s$. But then it follows by embedding that this is also an approximation in $F^s_{1,\infty} (\rn)$ for any $s$. In
	particular one has by Step 2
	\begin{\eq}   \label{2.51}
		W^\alpha_t: \quad \os{\circ}{F}{}^s_{1,q} (\rn) \hra \os{\circ}{F}{}^{s+d}_{1,q} (\rn), \qquad 1\le q \le \infty, \quad s\in \real,
		\quad d>0.
	\end{\eq}
	The operator $W^\alpha_t$  is self--dual, $(W^\alpha_t )' = W^\alpha_t$. Then \eqref{2.37} with $\As (\rn) = F^s_{\infty,q} (\rn)$,
	$s\in \real$, $1\le q <\infty$, follows from \eqref{2.50}, \eqref{2.51}.
\end{proof}
\begin{remark}\label{lift}
	We justify \eqref{2.49}. Let $f\in\Aspq,~1\le p<\infty,~s>0$ and $\sigma>0$. Without loss of generality we may assume $\sigma<\infty$ (otherwise one replaces $s$ by $s-\varepsilon$ with $0<\varepsilon<s$). If $f\in S(\rn)$ then
	\begin{\eq}  
		\Fou(W_t^{\alpha}I_\sigma f)(\xi)=e^{-t\lvert\xi\lvert^{2\alpha}}(1+\lvert\xi\lvert^2)^{-\sigma/2}(\Fou f)(\xi)
	\end{\eq}
	is well defined pointwise and belongs to $S'(\rn)$. Hence,
	\begin{\eq} 
		I_{-\sigma}	(W_t^{\alpha}I_\sigma f)=\Fou^{-1}e^{-t\lvert\xi\lvert^{2\alpha}}\Fou f=W_t^{\alpha} f.
	\end{\eq}	
	for all $f\in~S(\rn)$. Then it follows \eqref{2.49} for all $f\in\Aspq$ by \eqref{2.37} (with $d=0$), the lift property \eqref{2.37} and density of $S(\rn)$ in $\Aspq$.
\end{remark}	

\begin{remark}\label{temper}
	We observe that $\|f\lvert\,\accentset{\ast}{A}^s_{p,q}(\rn)\| $ is an equivalent norm in $\,\dot{A}^s_{p,q}(\rn)$  if $2\alpha k>s$ for Besov spaces and $2\alpha k>s+n$ for Triebel-Lizorkin spaces. This is a direct consequence of Remark \ref{R2.2} and Theorem \ref{T2.4}.
\end{remark}

\section{Nonlinear fractional heat equations}\label{S3}
In \cite{BaS19} we dealt with the Cauchy problem \eqref{1.4}, \eqref{1.5} where $\alpha>0$ is a natural number. The case $\alpha=1$ corresponds to a classical non-linear heat equation. We established mild and strong solutions in appropriate function spaces $\Lv\cap C^{\infty}(\rn\times(0,T))$ being fixed points of the operator $T_{u_0}$ defined in \eqref{1.7}.

The aim of this section is to extend some of these results to the case of fractional powers $\alpha >1/2$. This restriction results from the mapping properties of the non-linearity $Du^2$ in $\As$--spaces, see Proposition \ref{estinhom} below. In particular, we make use of the smoothing properties formulated in Theorem \ref{T2.6}.
\\

For later purposes we recall multiplication properties in the respective spaces $\Aspq$ derived in \cite{BaS19} including $p=\infty$ for $F$-spaces.
\begin{proposition}\label{estspace}
	Let $1\leq p,q\leq\infty$ and $(\frac{n}{p}-\frac{n}{2})_+<s<\infty$. Let $A\in\{B,\,F\}$. Then it holds
	\begin{\eq}\label{e-2.13}
		\|f\cdot g \lvert A^{s-(\np -s)_+-\varepsilon}_{p,q}(\rn)\|\leq c\,\|f\lvert\Aspq\|\cdot\|g\lvert\Aspq\|
	\end{\eq}
	for all $f,g\in\Aspq$ and all $\varepsilon >0$.       
\end{proposition}
\begin{proof}
	Spaces $\Aspq$ ($p<\infty$ for $F$-spaces) are multiplication algebras if $s>\frac{n}{p}$. For $F^s_{\infty,q}(\rn)$, $s>0$ this property follows from \cite[Thm. 2.41]{T20}. Then \eqref{e-2.13} holds with $\varepsilon\geq 0$. If $s=\np$ assertion \eqref{e-2.13} holds due to  \cite[Prop. 2.3]{BaS19} and embedding (2.24) in \cite{T20}. 
	Finally, if $s<\np$ then \eqref{e-2.13} is a consequence of \cite[Prop. 2.1]{BaS19} as well as the Sobolev-type embeddings in \cite[Prop. 2.2]{BaS19}.
\end{proof}

\begin{corollary}\label{corest}
	Let $1\leq p,q\leq\infty$, $(\frac{n}{p}-\frac{n}{2})_+<s<\infty$ and $\sigma<s-1-(\frac{n}{p}-s)_+$. Let $A\in\{B,\,F\}$. Then  it holds
	\begin{eqnarray}\label{estD}
	\|D(f\cdot g)\lvert A^\sigma_{p,q}(\rn)\|\leq c\,\|f\lvert \Aspq\|\cdot \|g\lvert\Aspq\|
	\end{eqnarray}
	for all $f,g\in\Aspq$.
\end{corollary}
\begin{proof}
	Clearly, we have 
	\[
	\|D(f\cdot g)\lvert A^\sigma_{p,q}(\rn)\|\leq c\,\|f\cdot g\lvert A^{\sigma+1}_{p,q}(\rn)\|.
	\]
	Thus, \eqref{estD} follows from Proposition \ref{estspace} and $\sigma +1 < s-\big(\np -s\big)_+$.
\end{proof}	
Next we derive an estimate of $\Tu $ as defined in \eqref{1.7} for fixed $t>0$ in appropriate function spaces $\Aspq$.	

\begin{proposition}\label{estinhom}
	Let $2\leq n\in\nat$, $1\leq p\leq\infty,~1\leq q\leq\infty$,   $\left(\frac{n}{p}-\frac{n}{2}\right)_+<s<\infty$ and let $\alpha>1/2$. Let $T>0$ and
	let $a,~v,~d$ such that
	\begin{\eq}\label{Bed}
		\frac{1}{\alpha}< v\leq\infty,\quad-\infty<a + \frac{1}{v}<\alpha ,\quad 1+\left(\frac{n}{p}-s\right)_+ < d<2\left(\alpha-\frac{1}{v}\right).
	\end{\eq}
	If  
	\begin{\eq}
		u_0\in\Ao~~\mbox{with}~~s_0\leq s~~\mbox{and}~~ u\in\Lav{2\alpha v}{\frac{a}{2\alpha}}{\Aspq}
	\end{\eq}
	then there exists a constant $c>0$, independent of $u_0$ and $u$, such that 
	\begin{eqnarray}\label{Tuo}
	\|\Tu u(\cdot, t)\lvert\Aspq\|&\leq & c\,t^{-\frac{s-s_{0}}{2\alpha}} \|u_0\lvert\Ao\|\\
	&+&c\,t^{1-\frac{1}{\alpha v}-\frac{d}{2\alpha}-\frac{a}{\alpha}}
	\|u\lvert\Lav{2\alpha v}{\frac{a}{2\alpha}}{\Aspq}\|^2.\notag
	\end{eqnarray}
	for all $t$ with $0<t<T$ (with $\frac{1}{v}=0$ and the modification \eqref{Lvinf} if $v=\infty$).
\end{proposition}
\begin{proof}
	Note that condition \eqref{Bed} with respect to $d$ implies $\alpha >\frac{1}{2}$.
	Using Theorem \ref{T2.6} and Corollary \ref{corest} with $s-d$ in place of $\sigma$ we can estimate as follows
	\begin{\eq}  
		\begin{aligned}
			&\|\Tu u(\cdot, t)\lvert\Aspq\|\\&\leq \|W_t^\alpha u_0\lvert\Aspq\|
			+	\int\limits_0^t\|W^\alpha_{t-\tau}Du^2(\cdot,\tau)\lvert\Aspq\|\,\text{d}\tau\\
			&\lesssim t^{-\frac{s-s_0}{2\alpha}}\|u_0\lvert\Ao\|
			+	\int\limits_0^t (t-\tau)^{-\frac{d}{2\alpha}}\|Du^2(\cdot,\tau)\lvert A^{s-d}_{p,q}(\rn)\|\,\text{d}\tau\label{68}\\
			&\lesssim 	t^{-\frac{s-s_0}{2\alpha}}\|u_0\lvert\Ao\|
			+	\int\limits_0^t (t-\tau)^{-\frac{d}{2\alpha}}\|u(\cdot,\tau)\lvert\Aspq\|^2\,\text{d}\tau.\\
		\end{aligned}
	\end{\eq}
	Here we used that
	\[
	\sigma = s-d<s-1-\big(\frac{n}{p}-s\big)_+
	\]
	according to \eqref{Bed}.
	By means of H\"{o}lder's inequality with exponent $\alpha v>1$ we obtain
	\begin{\eq}  \label{3.8}
		\begin{aligned}
			\int\limits_0^t (t-\tau)^{-\frac{d}{2\alpha}}\tau^{-\frac{a}{\alpha}}\tau^{+\frac{a}{\alpha}}
			&\,\|u(\cdot,\tau)\lvert\Aspq\|^2\,\text{d}\tau\\
			&\lesssim \,
			t^{1-\frac{1}{\alpha v}-\frac{d}{2\alpha}-\frac{a}{\alpha}}
			\|u\lvert\Lav{2\alpha v}{\frac{a}{2\alpha}}{\Aspq}\|^2.
		\end{aligned}
	\end{\eq}
	Here we used the conditions $a+\frac{1}{v}<\alpha$ as well as $d< 2(\alpha - \frac{1}{v})$ to ensure that the integral is finite. 
\end{proof}

\begin{theorem}\label{fix2}
	Let $2\leq n\in\nat, ~\frac{1}{2}<\alpha<\infty,~1\leq p\leq\infty,~1\leq q\leq\infty$  and let $A\in\{B,\,F\}$. Let
	\begin{\eq}\label{e-3.13}
		\np -2\alpha +1<s_0
	\end{\eq}
	and
	\begin{\eq} \label{e-3.14}
		\nph <s_0 + \alpha.
	\end{\eq}
	Let $u_0\in\Ao$. 
	\cm
	{\em (i)}
	Then there exists a number $T>0$ such that the Cauchy problem \eqref{1.4},\eqref{1.5}
	has a unique mild solution $u$ belonging to 
	\begin{eqnarray*}
		\Lav{2\alpha v}{\tfrac{a}{2\alpha}}{\Aspq},
	\end{eqnarray*}
	for all $s$ satisfying
	\begin{\eq}
		s_0 \leq s< s_0+ \min (\alpha,\,2\alpha -1) \label{e-3.15}
	\end{\eq}
	and
	\begin{\eq}\label{e-3.15a}
		s>\max(\nph,\, \big(\np-2\alpha+1\big)_+), 
	\end{\eq}
	where $a,~v$ such that
	\begin{\eq}\label{e-3.16}
		0\leq \frac{1}{v}<\frac{1}{2}\Big( 2\alpha -1 - \Big(\np -s\Big)_+\Big)
	\end{\eq}
	and 
	\begin{\eq}\label{e-3.17}
		s-s_0 < a + \frac{1}{v}<\min{\Big(\alpha,\, 2\alpha -1-\Big(\frac{n}{p}-s\Big)_+\Big)}.
	\end{\eq}
	{\em (ii)}
	The mild solution $u$ obtained in part (i) also belongs to the space $L_\infty((0,T),\Ao)$.
	Moreover, if, in addition, $\max (p,q)<\infty$ 
	then the above solution 
	$u(\cdot)$ converges to $u_0$ with respect to the norm in $\Ao$ if $t\to 0+$ .
\end{theorem}
\begin{proof}
	First we observe that assumptions \eqref{e-3.13} and \eqref{e-3.15} imply that
	\[
	0\leq s-s_0 < \min{\Big(\alpha,\, 2\alpha -1-\Big(\frac{n}{p}-s\Big)_+\Big)}.
	\]
	Hence, conditions \eqref{e-3.16} and \eqref{e-3.17} make sense.\\
	{\em Step 1.}
	We  choose $d$ such that 
	\[
	1 + \Big(\frac{n}{p}-s\Big)_+ < d < \min\Big( 2\alpha - (a+\frac{1}{v}),\, 2(\alpha - \frac{1}{v}) \Big)
	\]
	which is possible due to conditions \eqref{e-3.16} and \eqref{e-3.17}. Then it follows from Proposition \ref{estinhom} that
	\begin{align*}
	\|\Tu u(\cdot, t)\lvert\Aspq\|&\leq  c\,t^{-\frac{s-s_{0}}{2\alpha}} \|u_0\lvert\Ao\|\\
	&+c\,t^{1-\frac{1}{\alpha v}-\frac{d}{2\alpha}-\frac{a}{\alpha}}
	\|u\lvert\Lav{2\alpha v}{\tfrac{a}{2\alpha}}{\Aspq}\|^2.\notag
	\end{align*}
	for $0<t<T$.
	We multiply both sides with $t^{\frac{a}{2\alpha}}$. Raising to the power of $2\alpha v$ and integrating over $(0,T)$ yield
	\begin{align}\label{s11}
	&\int_0^T t^{av}\|\Tu u(\cdot,t)\lvert A^s_{p,q}(\rn)\|^{2\alpha v}\text{d}t\\
	&\leq  c\,T^{\delta}\,\|u_0\lvert \Ao\|^{2\alpha v}\nonumber
	+T^{\varkappa}\,\|u\lvert\Lav{2\alpha v}{\tfrac{a}{2\alpha}}{\Aspq}\|^{4\alpha v}
	\end{align}	
	with
	\begin{align}\label{s12a}
	\delta=(a-s+s_0)v+1>0,~~\mbox{since }~ s-s_0<a+\frac{1}{v}
	\end{align}
	and 
	\begin{align}\label{s12b}
	\varkappa=\Big(2\alpha -\frac{1}{v}-d-a\Big)v>0,~~\mbox{since } d < 2\alpha -(a+\frac{1}{v}).
	\end{align}
	Thus, $\Tu$ maps the unit ball $U_T$ in $\Lav{2\alpha v}{\frac{a}{2\alpha}}{\Aspq}$ into itself if $T$ is sufficiently small. \\
	As for the contraction property consider $u,v\in U_T$. A similar calculation with $d$ as above (cf. also \eqref{68} and \eqref{3.8}) yields
	\begin{align*}
	&\|\Tu u(\cdot,t)-\Tu v(\cdot,t)\lvert A^s_{p,q}(\rn)\|\\
	&\leq c\,t^{1-\frac{1}{\alpha v}-\frac{d}{2\alpha}-\frac{a}{\alpha}}
	\left(\int\limits_0^t \tau^{\alpha v}\|u^2(\cdot,\tau)-v^2(\cdot,\tau)\lvert A^{s-d+1}_{p,q}(\rn)\|\,\text{d}\tau\right)^{1/\alpha v}.
	\end{align*}
	Application of Proposition \ref{estspace} leads in combination with H\"{o}lder's inequality to
	\begin{align}\label{int}
	&\|\Tu u(\cdot,t)-\Tu v(\cdot,t)\lvert A^s_{p,q}(\rn)\|\leq  c\,t^{1-\frac{1}{\alpha v}-\frac{d}{2\alpha}-\frac{a}{\alpha}}\nonumber\\
	&\times\left(\int\limits_0^t\tau^{av}\|u(\cdot,\tau)-v(\cdot,\tau)\lvert \Aspq\|^{\alpha v}\text{d}\tau\right)^{1/2\alpha v}\\
	&\times\left(\int\limits_0^t\tau^{av}\|u(\cdot,\tau)+v(\cdot,\tau)\lvert \Aspq\|^{\alpha v}\,\text{d}\tau\right)^{1/2\alpha v}\nonumber
	\end{align}
	Let temporarily
	$X^s_T=\Lav{2\alpha v}{\tfrac{a}{2\alpha}}{\Aspq}$, then it follows from \eqref{int}
	\begin{align}
	\|\Tu u-\Tu v\lvert X^s_T\|\leq c\, T^{\frac{\vk}{2\alpha v}}\|u+v\lvert X^s_T\|\|u-v\lvert X^s_T\|\label{s50}
	\end{align}
	with the same $\varkappa$ as in \eqref{s12b}. 
	If $T>0$ is small enough, then $\Tu:U_T\mapsto U_T$ is a contraction. Since we deal with Banach spaces we have shown that $Tu$ has a unique fixed point in $U_T$ and hence a mild solution of the Cauchy problem \eqref{1.4},\eqref{1.5}.\\
	To extend the uniqueness to the whole space we proceed similarly to e.g. \cite{MG98}, \cite{Wu}. Let $u\in U_T$ be the above solution and $v\in X^s_T$ a second solution. We observe that \eqref{s50} holds 
	for any $0<t\leq T_0\leq T$. With $u\in U_T$ we obtain
	\begin{align}\label{s6}
	\|u-v\lvert X^s_{T_0}\|\leq c\, T_0^{\varkappa}\,(1+\|v\lvert X^s_T\|)\|u-v\lvert X^s_{T_0}\|.
	\end{align}
	If we choose $T_0>0$ small enough such that $c\, T_0^{\varkappa}\,(1+\|v\lvert X^s_T\|)<1$ it follows that $u(\cdot,t)=v (\cdot,t)$ for any $t\in(0,T_0]$. 
	Now we take $u(\cdot,T_0)\in \Aspq\hookrightarrow \Aso$ as new initial value and proceed as in the previous steps until \eqref{s6} inclusively.
	There exists a unique solution $\widetilde{u}$ in a neighbourhood $U_\delta(T_0)$ with $\widetilde{u}(\cdot,T_0)=u(\cdot,T_0)$. Since it holds that $\widetilde{u}(\cdot,t)=u(\cdot,t)$ for all $t\in(0,T_0]\cap U_\delta(T_0)$ we have extended $u$ to some interval $(0,T_1]$ with $T_0<T_1$. 
	Thus, we have prolongated $u(\cdot,t)-v(\cdot,t)=0$ to some interval $(0,T_1]$ where $T_0<T_1\leq T$ By iteration it follows the uniqueness in  $\Lav{2\alpha v}{\frac{a}{2\alpha}}{\Aspq}$.\\ 
	{\em Step 2.}
	We show part (ii) of the theorem. 
	To this end we first prove that the mild solution obtained in Step 1 belongs to $L_\infty((0,T), \Ao)$.
	
	Let $u_0\in\Ao$ and let $u\in L_\infty \big((0,T),\,\frac{a}{2\alpha},\,\Aspq\big)$ be the corresponding solution, where 
	$s$ and $a$ satisfy \eqref{e-3.15}, \eqref{e-3.15a} and \eqref{e-3.17}, where $\frac{1}{v}$ is replaced by 0. Let $0<t<T$. 
	It holds
	\begin{align*}
	\|u(\cdot,t)\lvert\Ao\|\notag
	\leq  \|\Sm{t}u_0\lvert\Ao\|+\int_0^t\|\Sm{t-\tau}Du^2(\cdot,\tau)\lvert\Ao\|\text{d}\tau.
	\end{align*}
	Taking into account \eqref{2.49} and the lift property \eqref{2.47} we may assume $s_0>0$.
	Concerning the first summand we obtain
	\begin{align}
	&  \|\Sm{t}u_0\lvert\Ao\|\notag\\
	\lesssim &  \|\int\limits_{\rn}\left(e^{-t\lvert\xi\lvert^{2\alpha}}\right)^\vee(x-y)\,u_0(y)\,\text{d}y\lvert \Ao\|\notag\\
	= & \|\int\limits_{\rn} \left(e^{-\lvert\xi\lvert^{2\alpha}}\right)^\vee(z)\,u_0(x-t^{1/2\alpha}z)\,\text{d}z\lvert\Ao\|\notag \\
	\leq & \int\limits_{\rn}\left(e^{-\lvert\xi\lvert^{2\alpha}}\right)^\vee(z)\|u_0(x-t^{1/2\alpha}z)\lvert\Ao\|\,\text{d}z\label{mink}\\
	\lesssim\, & \|u_0\lvert\Ao\| \int\limits_{\rn}\left(e^{-\lvert\xi\lvert^{2\alpha}}\right)^\vee(z)\,\text{d}z\label{trans}<\infty
	\end{align}
	independent of $t$, where \eqref{mink} follows from the generalized Minkowski inequality for Banach spaces and \eqref{trans} from the translation invariance of $\Aspq$-- spaces (see also Remark \ref{R2.2}) and \eqref{2.24}.\\
	In order to estimate the second summand we first consider the case that $s-s_0 \leq 1+\big(\np -s\big)_+$ and we put
	\[
	d:= 1+\big(\np -s\big)_+ -(s-s_0) + \ve,
	\]
	where $\ve$ is chosen such that $0<\ve < 2\alpha -1 - (\np -s\big)_+ - (s-s_0)$ according to \eqref{e-3.17}.
	Then, $d>0,~s-s_0 <\alpha -\frac{d}{2}$, and we may choose $a$ such that
	\[
	0< s-s_0 < a< \min \Big( \alpha-\frac{d}{2},\,2\alpha -1 - \big(\np -s\big)_+\Big).
	\]
	Applying Theorem \ref{T2.6} with $s_0-d$ in place of $s$ and Corollary \ref{corest} with $\sigma = s_0-d < s-1 - \big(\np -s\big)_+$ we obtain
	\begin{align}
	\int_0^t&\|\Sm{t-\tau} Du^2(\cdot,\tau)\lvert\Ao\|\text{d}\tau\notag\\
	& \leq c\,\int_0^t (t-\tau)^{-\frac{d}{2\alpha}}\|Du^2(\cdot,\tau)\lvert A^{s_0-d}_{p,q}(\rn)\|\text{d}\tau\notag\\
	&\leq c\,\int_0^t (t-\tau)^{-\frac{d}{2\alpha}}\|u(\cdot,\tau)\lvert\Aspq\|^2\text{d}\tau\notag\\
	&\leq c\, t^{-\frac{d}{2\alpha}-\frac{a}{\alpha}+1}
	\Big(\sup_{0<\tau<T} \tau^{\frac{a}{2\alpha}}\,\|u(\cdot,\tau)\lvert\Aspq\|\Big)^2 \label{linf}
	\end{align}
	because of $0<d<2\alpha$ and $a<\alpha $. 
	If $s-s_0 > 1+\big(\np -s\big)_+$ then we can apply Theorem \ref{T2.6} with $d=0$ and Corollary \ref{corest} with $s_0$
	in place of $\sigma$ to get
	\begin{align}\label{3.17a}
	&\int_0^t\|\Sm{t-\tau} Du^2(\cdot,\tau)\lvert\Ao\|\text{d}\tau \leq c\sup\limits_{0<t<T} \int_0^t \|u(\cdot,\tau)\lvert\Aspq\|^2\,\text{d}\tau \notag\\
	&~~~~~~~~~\leq c\, t^{-\frac{a}{\alpha}+1}\,\Big(\sup_{0<\tau<T} \tau^{\frac{a}{2\alpha}}\,\|u(\cdot,\tau)\lvert\Aspq\|\Big)^2~.
	\end{align}
	The boundedness of $u(t,\cdot)$ on $(0,T)$ in $\Ao$ follows from \eqref{trans}, \eqref{linf} (due to $a<\alpha-\frac{d}{2}$), and \eqref{3.17a}
	(because of $a<\alpha$).	
	Next we consider the limit of $u(t,\cdot)$ in $\Ao$ if $t\to 0+$. It holds 
	\begin{align*}
	& \|u(\cdot,t)-u_0\lvert\Ao\|\notag\\
	&\leq \, \|\Sm{t}u_0-u_0\lvert\Ao\|+\int_0^t \|\Sm{t-\tau}Du^2(\cdot,\tau)\lvert\Ao\|\text{d}\tau.
	\end{align*}
	The second summand on the right-hand side tends to zero if $t\to 0+$ as a consequence of  \eqref{linf}, \eqref{3.17a} and the conditions with respect to $a,~\alpha,$ and $d$.
	Using the identity 
	\begin{align*}
	u_0(x)= u_0(x)\cdot\Big(\Big( e^{-t\lvert\xi\lvert^{2\alpha}}\big)^\vee\Big)^\wedge(0)=(2\pi)^{-n/2}\,\int\limits_{\rn} \left(e^{-\lvert t\xi\lvert^{2\alpha}}\right)^\vee (x-y) u_0(x) \,\text{d}y
	\end{align*}
	we obtain the estimate
	\begin{align}
	&\|\Sm{t}u_0-u_0\lvert\Ao\|\notag\\
	\lesssim & \int\limits_{\vert z\lvert>N }\left(e^{-\lvert\xi\lvert^{2\alpha}}\right)^\vee(z)\|u_0(x-t^{1/2\alpha}z)-u_0(x)\lvert\Ao\|\,\text{d}z\\
	+  &  \int\limits_{\lvert z \lvert\leq N }\left(e^{-\lvert\xi\lvert^{2\alpha}}\right)^\vee(z)\|u_0(x-t^{1/2\alpha}z)-u_0(x)\lvert\Ao\|\,\text{d}z.\label{sec}
	\end{align}
	The first summand is lower than $\varepsilon$ if we choose $N$ large enough. Fixing this $N$ the second summand tends to zero for $t$ tending to zero. This follows from the fact that the Schwartz space $S(\rn)$ is dense in $\Aspq$ if $\max(p,q)<\infty$ and the continuity of the translation (See also \cite[Subsection 1.2.d]{Gr04} for more details with respect to approximate identities). This completes the proof.
\end{proof}

In addition to the results of the previous part one may ask for well-posedness of the Cauchy problem. The notation well-posedness is not totally fixed in the literature
(see the comments in \cite[Subsection 6.2.5]{T13}) We adapt the standard notation, see e.g. \cite{BP}. The Cauchy problem is called locally well-posed if there exists a 
unique mild and strong solution according to Theorem \ref{fix2}. In addition it is required continuous dependence of the solutions with respect to initial data.
This means that for solutions $u_1$ and $u_2$ of \eqref{1.4}, \eqref{1.5} according to Theorem \ref{fix2} 
with respect to initial data $u_0^1$ and $u_0^2$, respectively, 
for any $\varepsilon>0$ there exists a $\delta>0$ and a time $T>0$ such that for all $0<t<T$ 
\begin{\eq}\label{eps}
	\|u_1(\cdot,t)-u_2(\cdot,t)\lvert\Ao\|\leq \varepsilon
\end{\eq}
holds if
\begin{\eq}\label{del}
	\|u^1_0-u^2_0\lvert\Ao\|\leq \delta~.
\end{\eq}
It is sufficient to consider solutions
$u\in\Lav{\infty}{\frac{a}{2\alpha}}{\Aspq}$. Recall that by construction of the solution as a fixed point of $\Tu$ we may assume $\|u\lvert\Lav{\infty}{\frac{a}{2\alpha}}{\Aspq}\|\leq 1$.\\

\par
\begin{theorem}\label{stablesol}
	Let $u_i\in L_\infty ((0,T_i),\,\frac{a}{2\alpha},\,\Aspq )$ ($i=1,2$) be solutions of \eqref{1.4}, \eqref{1.5} obtained in Theorem \ref{fix2} 
	with initial data $u^i_0\in\Ao$ in the corresponding time interval $(0,T_i)$. Let $\max (p,q)<\infty$.
	Then under the conditions of Theorem \ref{fix2} the Cauchy problem \eqref{1.4}, \eqref{1.5} is locally well-posed.
\end{theorem}
\begin{proof}
	Let $u_1$, $u_2$ be two solutions of \eqref{1.4}, \eqref{1.5} with corresponding initial data $u_0^1$, $u_0^2$. 
	We have
	\begin{align}
	& \|u_1(\cdot,t)-u_2(\cdot,t)\lvert\Ao\| \notag\\
	&\leq  \|\Sm{t}(u^1_0-u^2_0)\lvert\Ao\|
	+\int_0^t\|\Sm{t-\tau}\,D(u_1^2-u_2^2)(\cdot,\tau)\lvert\Ao\|\text{d}\tau\notag\\
	\end{align}
	To estimate the first summand of the right-hand side we use again Theorem \ref{T2.6} with $d=0$. 
	The second summand can be treated in the same way as in Step 2 of the proof of Theorem \ref{fix2} with $u^1 - u^2$ in place of $u$. Note, that $u^1\in L_\infty ((0,T_1),\,\frac{a}{2\alpha},\,\Aspq )$ and $u^2\in L_\infty ((0,T_2),\,\frac{a}{2\alpha},\,\Aspq )$. Hence, application of Minkowski's inequality leads to
	\begin{align}\label{T}
	\|u_1(\cdot,t)-u_2(\cdot,t)\lvert\Ao\| 
	\leq c\,\|u_0^1 - u_0^2\lvert\Ao\|+c\,t^{-\frac{d}{2\alpha}-\frac{a}{\alpha}+1}
	\end{align}
	with the same choice of $d\geq 0$ as in step 2 of Theorem \ref{fix2} and
	for all $0<t<T$ with $ T\leq \min (T_1,T_2)$. Then the right hand side in \eqref{T} is lower the the given $\varepsilon$ if this $T$ is chosen small enough.
\end{proof}

\section{Comments and special cases}\label{S4}
As already mentioned in the introduction our approach allows to deal with the Cauchy problem \eqref{1.4}, \eqref{1.5} for initial data $u_0$ belonging to spaces $\Ao$ with smoothness $s_0$
satisfying the a-priori condition \eqref{e-3.13}, i.~e. $s_0>\np -2\alpha +1$. \\
It refers to the so-called supercritical case, where the existence of local (small $T$) solutions can be expected. For a detailed discussion of (sub/super)critical spaces in the context of Navier-Stokes and the related scalar nonlinear heat equation \eqref{1.4} we refer to \cite[Subsection 5.5]{T14} ($\alpha =1$), \cite[Subsection 3.2]{BaS19} ($\alpha\in\nat$) and the references given there.
The arguments can be adapted to the case of fractional $\alpha >1/2$.\\
The second a-priori condition \eqref{e-3.14}, i.~e. $\nph <s_0+\alpha$ is due to the mapping properties of the nonlinearity $Du^2$ in \eqref{1.4} and relevant if $\alpha >1$.
A breaking point is $\alpha = \nh +1$. If $1<\alpha \leq \nh +1$ then we obtain a dependence on the parameter $p$ whether or not all supercritical spaces are admitted. If $\alpha >\nh +1$ then the supercritical case can never be completely covered by our method.\\
A further notable exponent is  $\alpha = \frac{n+2}{4}~~(\alpha =\frac{5}{4}~\mbox{if}~n=3)$. In this case we have the coincidence $\np-2\alpha +1=\np - \nh$  resulting in consequences with respect to
spatial smoothness $s$ of our solution spaces $\Lav{2\alpha v}{\tfrac{a}{2\alpha}}{\Aspq}$. All these aspects will be discussed in more detail in the following.\\
We always assume $n\in\nat$, $n\geq 2$, $\alpha>\frac{1}{2}$, $1\leq p,q\leq\infty$, $A\in\{B,\,F \}$ and $a,~v$ as in \eqref{e-3.16} and \eqref{e-3.17}. In the figures below the area of admitted $s_0$ is shaded, that one of $s$ is hatched.
\begin{remark}\label{rem-4.1}
	The case $\frac{1}{2}<\alpha\leq1$. \\
	Let $s_0>\frac{n}{p}-2\alpha+1$ and let $u_0\in\Ao$. 
	Then \eqref{e-3.14} is satisfied and \eqref{e-3.15} and \eqref{e-3.15a} read as
	\[
	s_0\leq s <s_0+ 2\alpha -1~~~\mbox{and}~~~ s>\Big(\np - 2\alpha +1 \Big)_+,
	\]
	respectively.
	\begin{itemize}
		\item[\rm{(i)}] If $\np - 2\alpha +1 < s_0\leq 0$ and $\frac{n}{2\alpha -1}< p \leq\infty$ then 
		there exists a unique mild solution $u\in\Lva$, 
		where $ 0<s< s_0 + 2\alpha-1$.
		\item[\rm{(ii)}] If $ s_0 > \Big(\np - 2\alpha +1 \Big)_+$ and $1\leq p\leq\infty$ 
		then there exists a unique mild solution $u\in\Lva$, 
		where $ s_0\leq s< s_0 + 2\alpha-1$.
	\end{itemize}
	This is well-known in the case $\alpha =1$ (see, for example, \cite{Ba15}, \cite{BaS17}, and \cite[Subsection 4.4]{T14}).
	Supercritical spaces $\Ao$ are completely covered for all $p,~1\leq p\leq\infty$ (see Figure 1).
	In particular, initial data $u_0\in A^0_{p,q}(\rn)$ are admitted if $\frac{n}{2\alpha -1} <p \leq\infty.$\\
\end{remark}
\begin{remark}\label{rem-4.2}
	The case  $1<\alpha \leq\frac{n+2}{4}$.\\
	This case implies $n>2$ and \eqref{e-3.13} as well as \eqref{e-3.14} are satisfied if 
	\[
	s_0>\Big(\np -\alpha + 1\Big)_+ -\alpha .
	\]
	Conditions \eqref{e-3.15} and \eqref{e-3.15a} read as
	\[
	s_0\leq s <s_0+ \alpha~~~\mbox{and}~~~ s>\big(\np - 2\alpha +1 \big)_+,
	\]
	respectively.
	\begin{itemize}
		\item[\rm{(i)}] If $\Big(\np -\alpha + 1\Big)_+ -\alpha < s_0\leq 0$ and $\frac{n}{2\alpha -1}< p \leq\infty$  
		then there exists a unique mild solution $u\in\Lva$, 
		where $ 0<s< s_0 + \alpha$.
		\item[\rm{(ii)}] If $ s_0 > \Big(\np - 2\alpha +1 \Big)_+$ and $1\leq p\leq\infty$ 
		then there exists a unique mild solution $u\in\Lva$, 
		where $ s_0\leq s< s_0 + \alpha$.
	\end{itemize}
	Supercritical spaces $\Ao$ are completely covered if $1\leq p\leq \frac{n}{\alpha -1}$ (see Figure 2).
	In particular, initial data $u_0\in A^0_{p,q}(\rn)$ are admitted if $\frac{n}{2\alpha -1} <p \leq\infty.$\\
\end{remark}

\begin{figure}[htb]
	\centering
	\begin{minipage}[t]{.5\linewidth}
		\caption{$\frac{1}{2}<\alpha\leq 1$}
		\centering
		\begin{tikzpicture}[%
		scale=3.5, samples=50 
		]
		\def \bsp1{(0,-2/3)  -- (1,4/3) -- (0,4/3) }
		\fill[green!20!white] \bsp1;
		\def \bsp2{(0,0) -- (1/3,0) -- (1,4/3) -- (0,4/3) }
		\pattern[pattern=my north east lines, pattern color=blue] \bsp2;
		\draw[->,thick] (0,0) -- (1.2,0) node[right] {\tiny $\frac{1}{p}$};
		\draw[->,thick] (0,-.9) -- (0,1.6) node[above] {\tiny $s$};
		\draw[red, thick] plot[id=f1,domain=0:1] (\x,{2*\x-2/3});
		\draw[blue, thick, dashed](0,0) -- (1/3,0);
		\draw[blue, thick](1/3,0) -- (1/2,0);
		\draw[blue, thick] plot[id=f1,domain=1/2:1] (\x,{2*\x-1});
		\draw[green, thick,  dashed] plot[id=f1,domain=1/6:1] (\x,{2*\x-2/3});
		\draw[semithick, dashed](1/3,-1/3) -- (1/3,4/3);
		\node [label={\textcolor{blue}{\tiny $s>\left(\frac{n}{p}-2\alpha+1\right)_+$}}] at (1,1.3) {};
		\node [label={\textcolor{blue}{\tiny $s=\left(\frac{n}{p}-\frac{n}{2}\right)_{+}$}}] at (1.25,0.6) {};
		\node [label={\textcolor{PineGreen}{\tiny $s_0>\frac{n}{p}-2\alpha+1$}}] at (1.04,0.3) {};
		\node [label={\textcolor{red}{\tiny $s=\frac{n}{p}-2\alpha+1$}}] at (1/3+0.02,-2/3) {};
		\node [label={\textcolor{black}{\tiny $\frac{1}{p}=\frac{2\alpha-1}{n}$}}] at (1/6+0.075,1.322) {};
		\draw (1/2,0) -- (1/2,0) node[below=2pt] {\tiny $\frac{1}{2}$};
		\draw (1,-.02) -- (1,.02) node[below=2pt] {\tiny $1$};
		\draw [blue](0.4,1.0) -- (0.4,1.0) node[below=1pt] {\tiny $s$};
		\draw [PineGreen](1/6+.01,-1/6) -- (1/6+.01,-1/6) node[left=1pt] {\tiny $s_0$};
		\draw (-.02,-1/3) -- (.02,-1/3) node[left=2pt] {\tiny $s_0$};
		\draw (-.02,1/3) -- (.02,1/3) node[left=2pt] {\tiny $s$};
		\draw (-.02,-2/3) -- (.02,-2/3) node[left=2pt] {\tiny $-2\alpha+1$};
		\draw [decorate,decoration={brace,amplitude=10pt},xshift=0pt,yshift=0pt]
		(0,-1/3) -- (0,1/3) node [black,midway,xshift=-1cm] 
		{\tiny $2\alpha-1$};
		\end{tikzpicture}
	\end{minipage}%
	\hfill%
	\begin{minipage}[t]{.5\linewidth}
		\caption{$1<\alpha\leq\frac{n+2}{4}$}
		\centering	
		\begin{tikzpicture}[%
		scale=3.5, samples=50 
		]
		\def \bsp1{(0,-1/3) -- (1/6, -1/3) -- (1,4/3) -- (0,4/3) }
		\fill[green!20!white] \bsp1;
		\def \bsp2{(0,0) -- (1/8,0) -- (19/24,4/3) -- (0,4/3) }
		\pattern[pattern=my north east lines, pattern color=blue] \bsp2;
		\draw[->,thick] (0,0) -- (1.2,0) node[right] {\tiny $\frac{1}{p}$};
		\draw[->,thick] (0,-.9) -- (0,1.6) node[above] {\tiny $s$};
		\draw[PineGreen, thick, dashed](0,-1/3) -- (1/6,-1/3);
		\draw[red, thick] plot[id=f1,domain=0:1] (\x,{2*\x-2/3});
		\draw[blue, thick] plot[id=f1,domain=1/8:19/24] (\x,{2*\x-1/4});
		\draw[blue, thick] plot[id=f1,domain=0:1/8] (\x,{0});
		\draw[blue, thick, dashed](0,0) -- (1/3,0);
		\draw[blue, thick](1/3,0) -- (1/2,0);
		\draw[blue, thick] plot[id=f1,domain=1/2:1] (\x,{2*\x-1});
		\draw[PineGreen, thick,  dashed] plot[id=f1,domain=1/6:1] (\x,{2*\x-2/3});
		\draw[semithick, dashed](1/8,-1/3) -- (1/8,4/3);
		\draw[semithick, dashed](1/3,0) -- (1/3,4/3);
		\node [label={\textcolor{blue}{\tiny $s>\left(\frac{n}{p}-\alpha+1\right)_+$}}] at (1.1,1.3) {};
		\node [label={\textcolor{PineGreen}{\tiny $s_0> -\alpha+\left(\frac{n}{p}-\alpha+1\right)_+$}}] at (1.4,1) {};
		\node [label={\textcolor{blue}{\tiny $s=\left(\frac{n}{p}-\frac{n}{2}\right)_{+}$}}] at (1.25,0.6) {};
		\node [label={\textcolor{red}{\tiny $s=\frac{n}{p}-2\alpha+1$}}] at (1/3+0.02,-2/3) {};
		\node [label={\textcolor{black}{\tiny $\frac{1}{p}=\frac{\alpha-1}{n}$}}] at (1/6+0.075,1.322) {};
		\draw (1/2,0) -- (1/2,0) node[below=2pt] {\tiny $\frac{1}{2}$};
		\draw (1/3,0) -- (1/3,0) node[below=2pt] {\tiny $\frac{2\alpha-1}{n}$};
		\draw (1,-.02) -- (1,.02) node[below=2pt] {\tiny $1$};
		\draw [blue](0.4,1.0) -- (0.4,1.0) node[below=1pt] {\tiny $s$};
		\draw [PineGreen](1/6+.01,-1/6) -- (1/6+.01,-1/6) node[left=1pt] {\tiny $s_0$};
		\draw (-.02,1/6) -- (.02,1/6) node[left=2pt] {\tiny $s$};
		\draw (-.02,-1/6) -- (.02,-1/6) node[left=2pt] {\tiny $s_0$};
		\draw (-.02,-1/3) -- (.02,-1/3) node[left=2pt] {\tiny $-\alpha$};
		\draw (-.02,-2/3) -- (.02,-2/3) node[left=2pt] {\tiny $-2\alpha+1$};
		\draw [decorate,decoration={brace,amplitude=10pt},xshift=0pt,yshift=0pt]
		(0,-1/6) -- (0,1/6) node [black,midway,xshift=-.5cm] 
		{\tiny $\alpha$};
		\end{tikzpicture}
	\end{minipage}
\end{figure}

\newpage
\begin{remark}\label{rem-4.3}
	The case $\frac{n+2}{4}<\alpha \leq\nh +1.$\\
	As in the previous case 
	\eqref{e-3.13} as well as \eqref{e-3.14} are satisfied if 
	\[
	s_0>\Big(\np -\alpha + 1\Big)_+ -\alpha .
	\]
	Now, conditions \eqref{e-3.15} and \eqref{e-3.15a} read as
	\[
	s_0\leq s <s_0+ \alpha~~~\mbox{and}~~~ s>\nph,
	\]
	respectively.
	\begin{itemize}
		\item[\rm{(i)}] If $\Big(\np -\alpha + 1\Big)_+ -\alpha < s_0\leq \nph$ and $1\leq p \leq\infty$  
		then there exists a unique mild solution $u\in\Lva$, 
		where $ \nph <s< s_0 + \alpha$.
		\item[\rm{(ii)}] If $ s_0 > \Big(\np - 2\alpha +1 \Big)_+$ and $1\leq p\leq\infty$ 
		then there exists a unique mild solution $u\in\Lva$, 
		where $ s_0\leq s< s_0 + \alpha$.
	\end{itemize}
	Supercritical spaces $\Ao$ are completely covered if $1\leq p\leq \frac{n}{\alpha -1}$ (see Figure 3).
	Initial data $u_0\in A^0_{p,q}(\rn)$ are admitted provided that $\frac{n}{2\alpha -1} <p \leq\infty$ if $\frac{n}{2\alpha -1}\geq 1$
	or $1\leq p\leq\infty$ if $\frac{n}{2\alpha -1}<1$.\\
\end{remark}
\begin{remark}\label{rem-4.4}
	The case $\alpha >  \nh +1.$\\
	In this case 
	\eqref{e-3.14} implies \eqref{e-3.13}. We have  
	\[
	s_0>\nph  -\alpha > \np -2\alpha +1 .
	\]
	Conditions \eqref{e-3.15} and \eqref{e-3.15a} read as
	\[
	s_0\leq s <s_0+ \alpha~~~\mbox{and}~~~ s>\nph,
	\]
	respectively.
	\begin{itemize}
		\item[\rm{(i)}] If $\nph -\alpha < s_0\leq \nph$ and $1\leq p \leq\infty$  
		then there exists a unique mild solution $u\in\Lva$, 
		where $ \nph <s< s_0 + \alpha$.
		\item[\rm{(ii)}] If $ s_0 > \nph$ and $1\leq p\leq\infty$ 
		then there exists a unique mild solution $u\in\Lva$, 
		where $ s_0\leq s< s_0 + \alpha$.
	\end{itemize}
	Supercritical spaces $\Ao$ can never be  completely covered for given $p,~1\leq p\leq\infty$ (see Figure 4).
	Initial data $u_0\in A^0_{p,q}(\rn)$ are admitted for all $p,~1\leq p\leq\infty$.\\
\end{remark}

\begin{figure}[htb]
	\centering
	\begin{minipage}[t]{.5\linewidth}
		\caption{$\frac{n+2}{4}<\alpha\leq \frac{n}{2}+1$}
		\centering
		\begin{tikzpicture}[%
		scale=4, samples=50 
		]
		\def \bsp1{(0,-2/3) -- (1/6, -2/3) -- (2/3,1/3) -- (2/3,2/3) -- (0,2/3) }
		\fill[green!20!white] \bsp1;
		\def \bsp2{(0,0) -- (1/3,0) -- (2/3,2/3) -- (2/3,2/3) -- (0,2/3) }
		\pattern[pattern=my north east lines, pattern color=blue] \bsp2;
		\draw[->,thick] (0,0) -- (2/3+.2,0) node[right] {\tiny $\frac{1}{p}$};
		\draw[->,thick] (0,-1.08) -- (0,0.8) node[above] {\tiny $s$};
		\draw[PineGreen, thick, dashed](0,-2/3) -- (1/6,-2/3);
		\draw[red, thick] plot[id=f1,domain=0:2/3] (\x,{2*\x-1});
		\draw[blue, thick, dashed](0,0) -- (1/3,0);
		\draw[blue, thick, dashed] plot[id=f1,domain=1/3:2/3] (\x,{2*\x-2/3});
		\draw[PineGreen, thick,  dashed] plot[id=f1,domain=1/6:2/3] (\x,{2*\x-1});
		\draw[semithick, dashed](1/6,-2/3) -- (1/6,2/3);
		\node [label={\textcolor{blue}{\tiny $s>\left(\frac{n}{p}-\frac{n}{2}\right)_{+}$}}] at (1.05,2/3-.2) {};
		\node [label={\textcolor{PineGreen}{\tiny $s_0>\left(\frac{n}{p}-\alpha+1\right)_+-\alpha$}}] at (1.07,0.2) {};
		\node [label={\textcolor{red}{\tiny $s=\frac{n}{p}-2\alpha+1$}}] at (1/2+0.01,-2/3) {};
		\node [label={\textcolor{black}{\tiny $\frac{1}{p}=\frac{\alpha-1}{n}$}}] at (1/6+0.075,0.65) {};
		\draw (1/2,0) -- (1/2,0) node[below=2pt] {\tiny $\frac{2\alpha-1}{n}$};
		\draw (1/3,0) -- (1/3,0) node[below=2pt] {\tiny $\frac{1}{2}$};
		\draw (2/3,-.02) -- (2/3,.02) node[below=2pt] {\tiny $1$};
		\draw [blue](0.3,1/3) -- (0.3,1/3) node[below=1pt] {\tiny $s$};
		\draw [PineGreen](1/6+.01,-1/6) -- (1/6+.01,-1/6) node[left=1pt] {\tiny $s_0$};
		\draw (-.02,-2/3) -- (.02,-2/3) node[left=2pt] {\tiny $-\alpha$};
		\draw (-.02,-1) -- (.02,-1) node[left=2pt] {\tiny $-2\alpha+1$};
		\end{tikzpicture}
	\end{minipage}%
	\hfill%
	\begin{minipage}[t]{.5\linewidth}
		\caption{$\alpha>\frac{n}{2}+1$}
		\centering	
		\begin{tikzpicture}[%
		scale=3.5, samples=50 
		]				
		
		\def \bsp1{(0,-1/2) -- (1/6, -1/2) -- (3/8,-1/12) -- (3/8,1) -- (0,1) }
		\fill[green!20!white] \bsp1;
		\def \bsp2{(0,0) -- (1/6,0) --(3/8,5/12) -- (3/8,1) -- (0,1) }
		\pattern[pattern=my north east lines, pattern color=blue] \bsp2;
		\draw[->,thick] (0,0) -- (5/8,0) node[right] {\tiny $\frac{1}{p}$};
		\draw[->,thick] (0,-1.08) -- (0,1.1) node[above] {\tiny $s$};
		\draw[PineGreen, thick, dashed](0,-1/2) -- (1/6,-1/2);
		\draw[red, thick] plot[id=f1,domain=0:3/8] (\x,{2*\x-6/6});
		\draw[blue, thick, dashed](0,0) -- (1/6,0);
		\draw[blue, thick, dashed] plot[id=f1,domain=1/6:3/8] (\x,{2*\x-1/3});
		\draw[PineGreen, thick,  dashed] plot[id=f1,domain=1/6:3/8] (\x,{2*\x-5/6});
		\draw[semithick, dashed](3/8,-1/4) -- (3/8,1);
		\node [label={\textcolor{blue}{\tiny $s>\left(\frac{n}{p}-\frac{n}{2}\right)_{+}$}}] at (.75,.6) {};
		\node [label={\textcolor{PineGreen}{\tiny $s_0>\left(\frac{n}{p}-\frac{n}{2}\right)_{+}-\alpha$}}] at (.8,1/3) {};
		\node [label={\textcolor{red}{\tiny $s=\frac{n}{p}-2\alpha+1$}}] at (.8,-3/8) {};
		\draw (1/6,0) -- (1/6,0) node[below=2pt] {\tiny $\frac{1}{2}$};
		\draw (1/3+.1,0) -- (1/3+.1,0) node[below=2pt] {\tiny $1$};
		\draw [blue](1/6,2/3) -- (1/6,2/3) node[below=1pt] {\tiny $s$};
		\draw [PineGreen](1/6,-1/4) -- (1/6, -1/4) node[left=1pt] {\tiny $s_0$};
		\draw (-.02,-1/2) -- (.02,-1/2) node[left=2pt] {\tiny $-\alpha$};
		\draw (-.02,-1) -- (.02,-1) node[left=2pt] {\tiny $-2\alpha+1$};
		\end{tikzpicture}
	\end{minipage}%
\end{figure}
\newpage
\begin{remark}\label{rem-4.5}
	Some attention attracted Cauchy problems of type \eqref{1.4}, \eqref{1.5} for fractional power dissipative equations
	with initial data belonging to spaces $L_p(\rn),~1<p<\infty$ (see, for example, \cite{MYZ08}).
	Let us suppose that initial data $u_0 \in A^0_{p,q}(\rn)$, where $1<p<\infty,~1\leq q\leq\infty$. 
	In particular, this applies to $u_0\in L_p(\rn)= F^0_{p,2}(\rn)$. We follow our approach and recall (see the preceding remarks)
	that $u_0\in A^0_{p,q}(\rn)$ is admitted if $0<\frac{1}{p}<\frac{2\alpha -1}{n}
	~~(0<\frac{1}{p}\leq 1~~\mbox{if}~~\frac{2\alpha -1}{n}>1)$ for given $\alpha >\frac{1}{2}$ and $n\in\nat~(n\geq 2)$.
	Introducing a new parameter $\mu = a +\frac{1}{v}$ we can reformulate
	\[
	u \in \Lav{2\alpha v}{\tfrac{a}{2\alpha}}{\Aspq} 
	\]
	in Theorem \ref{fix2} as 
	\begin{\eq}\label{e-4.1}
		\int\limits_0^T  t^{\mu v} \|\,u(\cdot, t)\lvert \As(\rn)\|^{2\alpha v}\, \frac{\text{d} t}{t}~<\infty
	\end{\eq}
	(with modification $\sup\limits_{0<t<T} t^{\frac{\mu}{2\alpha}}\| \dots\| <\infty$ if $v=\infty$),
	where 
	\begin{align*}
	& \nph <s <\min (\alpha, 2\alpha -1),\\
	& \mu >s ~~~\mbox{and} ~~~0\leq \frac{2}{v} <\Big(2\alpha -1 - \big(\np -s\big)_+\Big)~.
	\end{align*}
	In the following we shall make use of the embeddings
	\begin{\eq}\label{e-4.2}
		A^{\lambda +\np}_{p,q}(\rn)~~\hookrightarrow ~~A^{\lambda}_{\infty , q}(\rn)
	\end{\eq}
	for $\lambda \geq 0$ and
	\begin{\eq}\label{e-4.3}
		A^{\np-\nr}_{p,q}(\rn)~~\hookrightarrow ~~A^{0}_{r, q}(\rn)
	\end{\eq}
	for $0<\frac{1}{r}<\frac{1}{p}$. We always assume $u_0\in A^{0}_{p,q}(\rn)$ and $u$ stands for the unique mild solution
	according to Theorem \ref{fix2}. We distinguish the following cases:\\[1ex]
	{\bf Case 1:} \quad Let $\frac{1}{2}< \alpha <1$ and let $0<\frac{1}{p}< \frac{2\alpha -1}{n}$. It follows from 
	\eqref{e-4.1} and \eqref{e-4.2} that
	\begin{\eq}\label{e-4.4}
		\int\limits_0^T  t^{\mu v} \|\,u(\cdot, t)\lvert A^{\lambda}_{\infty,q}(\rn)\|^{2\alpha v}\, \frac{\text{d} t}{t}~<\infty,
	\end{\eq}
	where $0\leq \lambda < 2\alpha - 1-\np,~~\mu >\lambda +\np,$ and $0\leq \frac{2}{v}< 2\alpha -1$.\\
	Combining \eqref{e-4.1} and \eqref{e-4.3} we get
	\begin{\eq}\label{e-4.5}
		\int\limits_0^T  t^{\mu v} \|\,u(\cdot, t)\lvert A^{0}_{r,q}(\rn)\|^{2\alpha v}\, \frac{\text{d} t}{t}~<\infty,
	\end{\eq}
	for $0<\frac{1}{r}<\frac{1}{p},~~\mu >\np -\nr,$ and $0\leq \frac{2}{v}< 2\alpha -1 -\nr$.\\[2ex] 
	{\bf Case 2:} \quad Let $\alpha \geq 1$ and let $0<\frac{1}{p}< \frac{\alpha}{n}~~(\frac{1}{p}\leq 1~\mbox{if}~\alpha >n)$. 
	According to \eqref{e-4.1} and \eqref{e-4.2} we find that \eqref{e-4.4} holds for 
	$0\leq \lambda < \alpha -\np,~~\mu >\lambda +\np,$ and $0\leq \frac{2}{v}< 2\alpha -1$. Moreover, it holds \eqref{e-4.5}
	for  
	$0<\frac{1}{r}<\min\big(\frac{1}{p},\,\frac{1}{2}\big),$\\
	$\mu >\np -\nr > \nph,$ and $0\leq \frac{2}{v}< 2\alpha -1 -\nr$.\\[2ex] 
	{\bf Case 3:} \quad Let $\alpha \geq 1$ and let $\frac{\alpha}{n}\leq \frac{1}{p} <\frac{2\alpha -1}{n}~~
	(\frac{1}{p}\leq 1~\mbox{if}~\alpha >n)$. Then $u$ satisfies \eqref{e-4.5} for 
	$0<\frac{1}{p}-\frac{\alpha}{n}<\frac{1}{r}<\min\big(\frac{1}{p},\,\frac{1}{2}\big),~~
	\mu >\np -\nr > \nph,$ and $0\leq \frac{2}{v}< 2\alpha -1 -\nr$.\\[2ex] 
	We may choose $\mu =\frac{1}{v}$ in \eqref{e-4.5} if 
	\[
	\np -\nr<\frac{1}{v}< \frac{1}{2} (2\alpha -1 -\nr),~~\frac{1}{p}<\frac{2\alpha -1}{n} <1,~~
	0<\frac{1}{r}<\min\big(\frac{1}{p},\,\frac{1}{2}\big).
	\]
	For example, this is the case if $\frac{1}{p}\leq \frac{2\alpha -1}{2n} \leq\frac{1}{2}$. Then 
	\[
	\np -\nr<\frac{1}{v}< \frac{1}{2} (2\alpha - 1) -\nr < \frac{1}{2}(2\alpha -1-\nr) ~~\mbox{and}~~
	\frac{1}{p}<\frac{2\alpha -1}{2n} < \frac{\alpha}{n}.
	\] 
	Thus, it follows from Case 2 and \eqref{e-4.5} that
	\begin{\eq}\label{e-4.6}
		\int\limits_0^T   \|\,u(\cdot, t)\lvert A^{0}_{r,q}(\rn)\|^{2\alpha v}\, \text{d} t~<\infty,
	\end{\eq}
	if $u_0\in A^{0}_{p,q}(\rn),~ 0<\frac{1}{r}<\frac{1}{p}<\frac{2\alpha -1}{2n}\leq\frac{1}{2},$ and 
	$\frac{1}{v}>\np -\nr$.
	Results of type \eqref{e-4.5} (in the case $v=\infty$) and \eqref{e-4.6}  can be found in \cite[Theorems 4.3 and 4.4]{MYZ08}.
	In a certain sense Theorem \ref{fix2} is an extension of their investigations (in the case b=d=1).
\end{remark}

\end{document}